\documentclass[a4paper,11pt]{article}
\usepackage{titlesec}
\usepackage{amsthm}
\usepackage{amssymb}
\usepackage{color}
\usepackage{mathtools}
\usepackage{geometry}

\usepackage{algorithm2e}

\usepackage{verbatim}
\numberwithin{equation}{section}

\theoremstyle{definition}

\newcommand{\bfs}[1]{{\boldsymbol #1}}

\usepackage{array}
\newcolumntype{C}[1]{>{\centering\arraybackslash}m{#1}}
\usepackage{graphics}
\usepackage{amssymb, latexsym, amsmath, amsfonts, epsfig}
\usepackage{bm}
\usepackage{graphicx}
\usepackage{color}
\usepackage{epstopdf}
\usepackage{comment}
\usepackage{soul}
\usepackage[normalem]{ulem}
\usepackage{float}
\usepackage{MnSymbol}

\usepackage{caption}
\usepackage{amsmath,amssymb}

\usepackage{lipsum}
\usepackage{amsfonts}
\usepackage{graphicx}
\usepackage{epstopdf}
\usepackage{algorithmic}
\ifpdf
  \DeclareGraphicsExtensions{.eps,.pdf,.png,.jpg}
\else
  \DeclareGraphicsExtensions{.eps}
\fi

\def\d{{\, \rm d}}

\usepackage{enumitem}
\setlist[enumerate]{leftmargin=.5in}
\setlist[itemize]{leftmargin=.5in}

\title{Superfloe Parameterization with Physics Constraints for Uncertainty Quantification of Sea Ice Floes}

\author{Nan Chen\thanks{Department of Mathematics, University of Wisconsin--Madison, Madison, WI 53706, USA
  (Nan.Chen@wisc.edu)}
  \and Quanling Deng\thanks{Corresponding author. Department of Mathematics, University of Wisconsin--Madison, Madison, WI 53706, USA (Quanling.Deng@math.wisc.edu); Current address: School of Computing, Australian National University, Canberra, ACT 2601, Australia (Quanling.Deng@anu.edu.au) }
  \and Samuel N. Stechmann\thanks{Department of Mathematics and Department of Atmospheric and Oceanic Sciences, University of Wisconsin--Madison, Madison, WI 53706, USA (Stechmann@wisc.edu)}}

\usepackage{amsopn}

\begin{document}

\maketitle

\begin{abstract}
The discrete element method (DEM) is providing a new modeling approach for describing sea ice dynamics. It exploits particle-based methods to characterize the physical quantities of each sea ice floe along its trajectory under Lagrangian coordinates. One major challenge in applying the DEM models is the heavy computational cost when the number of floes becomes large. In this paper, an efficient Lagrangian parameterization algorithm is developed, which aims at reducing the computational cost of simulating the DEM models while preserving the key features of the sea ice. The new parameterization takes advantage of a small number of artificial ice floes, named the superfloes, to effectively approximate a considerable number of the floes, where the parameterization scheme satisfies several important physics constraints. The physics constraints guarantee the superfloe parameterized system will have similar short-term dynamical behavior as the full system. These constraints also allow the superfloe parameterized system to accurately quantify the long-range uncertainty, especially the non-Gaussian statistical features, of the full system. In addition, the superfloe parameterization facilitates a systematic noise inflation strategy that significantly advances an ensemble-based data assimilation algorithm for recovering the unobserved ocean field underneath the sea ice. Such a new noise inflation method avoids ad hoc tunings as in many traditional algorithms and is computationally extremely efficient. Numerical experiments based on an idealized DEM model with multiscale features illustrate the success of the superfloe parameterization in quantifying the uncertainty and assimilating both the sea ice and the associated ocean field.
 \textbf{Mathematics Subjects Classification}:   65P99, 65M99, 76M28
\end{abstract}



\section{Introduction}

\begin{figure}[h!]
	\centering
	\includegraphics[width=3in]{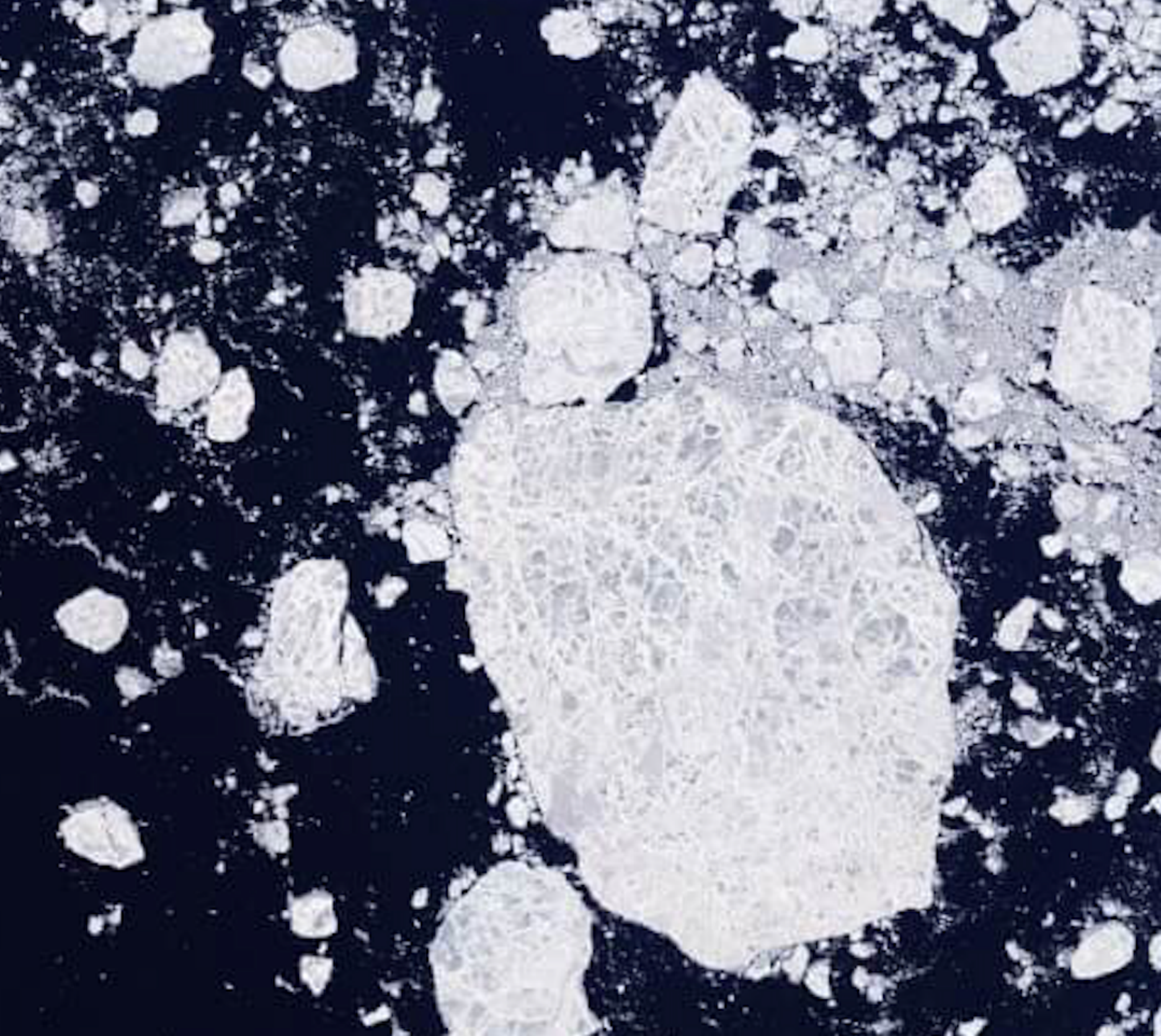}
	\caption{Arctic sea ice floes. NASA credits, figure from the Global Land Cover Facility.}
	\label{fig:sif}
\end{figure}

Sea ice forms when seawater freezes. It has a profound influence on the polar environment, such as weather, climate, and ocean circulation. It also interacts with many other climate phenomena across different time scales and therefore influences the entire earth system \cite{cavalieri2012arctic, comiso2008accelerated, parkinson2014global, serreze2007perspectives}. A sea ice floe is defined as a large pack of floating ice \cite{zhang2018sea}, which is widely observed in the marginal ice zone. See Figure \ref{fig:sif} for an example of sea ice floes in  the Arctic.

Sea ice is widely modeled as a continuum for its viscous-plastic rheology \cite{hibler1979dynamic, hunke1997elastic}, which is reasonable for describing the large-scale features of sea ice. However, at scales of the order of 10 km and smaller, sea ice exhibits brittle behavior with individual fragments clearly visible from satellite observations. To better characterize such granular media, discrete element method (DEM) models have become important in modeling sea ice, as an alternative to the continuum models. The DEM models exploit particle-based approaches to describe the physical quantities of each sea ice floe along its trajectory under the Lagrangian coordinates \cite{cundall1979discrete,hopkins2004discrete, bateman2019simulating, xu2012discrete,herman2016discrete, damsgaard2018application}. The DEM models have several advantages over the traditional continuum models \cite{lindsay2004new,tuhkuri2018review}.
First, the continuum models, which are based on Eulerian coordinates, often require a simulation domain that contains floe-free (open seawater) regions. This leads to a waste of computational resources. In contrast, the DEM models track the floes in the relevant regions, which significantly improves the computational efficiency and accuracy.
Second, the Lagrangian DEM models have a large flexibility in varying the spatial resolution, while the Eulerian models require adaptive meshes for this purpose, which often introduce additional computational costs \cite{zienkiewicz1991adaptivity}.
In addition, the observed trajectories of the floes can be directly incorporated into Lagrangian data assimilation (DA) to recover the underlying ocean flow field, which typically lacks direct observations in the polar regions.

Despite many computational merits over the traditional continuum models, computational efficiency remains a challenging issue of applying the DEM models when the number of the floes becomes large  \cite{herman2016discrete}. In addition to increased computational storage, frequent collisions between floes require a short numerical integration time step. The demanding computational cost also makes it extremely difficult to carry out efficient DA, which often requires running an ensemble of model forecasts. Therefore, developing effective parameterizations has become an important topic to facilitate computational efficiency and accuracy in applying DEM models.

In this paper, we develop an efficient physics-constrained superfloe parameterization scheme for the Lagrangian DEM models. It aims at utilizing  a small number of artificial ice floes, named superfloes, to effectively approximate a considerable number of small-scale floes. Notably, the new parameterization scheme satisfies several important physics constraints. The dimension of the resulting parameterized system is much lower than the full system and therefore the computational cost is significantly reduced.
The physics constraints guarantee the superfloe parameterized system to have similar short-term dynamical behavior as the full system. These constraints also allow the superfloe parameterized system to accurately quantify the long-range uncertainty, especially the non-Gaussian statistical features, of the full system.  The basic idea of such an efficient physics-constrained superfloe parameterization scheme is to iteratively group together the neighbouring small floes to create a superfloe that characterizes the key physics and statistics of the floe clusters. The superfloe parameterization shares some common features as the coarsening procedures in the smoothed particle hydrodynamics (SPH) method \cite{monaghan1992smoothed,monaghan2005smoothed,monaghan2012smoothed} in which reductions in computational cost can be obtained by merging the nearest small particles \cite{hu2017consistent,vacondio2013variable,barcarolo2014adaptive,vacondio2016variable}. Nevertheless, there are quite a few unique features of the superfloe parameterization and its application here.
First, the superfloes have clear physical meanings. Similar to the real floes, these superfloes are equipped with all the physical quantities, including the floe radius,
thickness, center of mass, velocity, and angular velocity. Second, the superfloes are designed to satisfy several key physics constraints. For example, the total mass and the sea ice concentration in the superfloe parameterized model are the same as those in the full model. The momentum and the angular momentum of each superfloe are also constrained to have the same value as those associated with the floe clusters. In addition, the statistics of all the kinematic quantities, as well as those of the floe-floe contact forces, are by design to be retained in the superfloe parameterized system. Third, the superfloe parameterization facilitates systematic noise inflation in Lagrangian DA \cite{apte2008bayesian, ide2002lagrangian, chen2014information}, which advances using the ensemble-based DA algorithm for recovering the ocean field underneath the sea ice. Such a new noise inflation method avoids ad hoc tunings as in the traditional approaches and is necessary for the ensemble-based DA algorithm to work effectively. Finally, the superfloe parameterization algorithm is simple to implement and is amenable to different DEM systems.

The rest of the paper is organized as follows. Section \ref{sec:ps} summarizes a recently developed idealized DEM floe model of sea ice. Section \ref{sec:sp} develops the superfloe parameterization algorithm. Section \ref{sec:uq} studies the short- and long-range uncertainty quantification of the superfloe parameterization.  Section \ref{sec:da} aims at showing the advantage of applying superfloe parameterization to facilitate the Lagrangian DA. The concluding remarks are included in Section \ref{sec:con}.

\section{An Idealized DEM Model for Sea-Ice Floe Dynamics} \label{sec:ps}
\subsection{Overview}
This section aims at summarizing a recently developed idealized DEM model for sea-ice floe dynamics \cite{chen2020}, which will be utilized as a test model for the uncertainty quantification (UQ) and DA of the superfloe parameterization.

The sea ice floes move in a two-dimensional square domain with double periodic boundary conditions for the ocean. Although the atmospheric forcing is the dominant contribution to the sea ice velocity, it often lies on the large scale. In other words, the main role of the atmospheric forcing is to provide a homogeneous background velocity for the floes at each fixed time instant, which from a mathematical point of view can be eliminated by applying a Galilean transformation. This means the large-scale atmospheric forcing has little impact on the superfloe parameterization. The atmospheric forcing is therefore omitted in the study here. Furthermore, in order to develop a simple and mathematically tractable DEM model for testing the superfloe parametrization skill, only the most fundamental dynamical features are included in the model. In particular, the model emphasizes the impact of sea ice motion from the ocean drag force as well as the collisions between the ice floes. The latter is a unique property that distinguishes sea ice floes from free tracers. On the other hand, some other sea ice characteristics, such as melting and ridging, are omitted in the current study for simplicity since they play a relatively less important role in affecting the superfloe parameterization. 

\subsection{The geometry of the floes}
\begin{figure}[h!]
	\centering
	\includegraphics[width=2.5in]{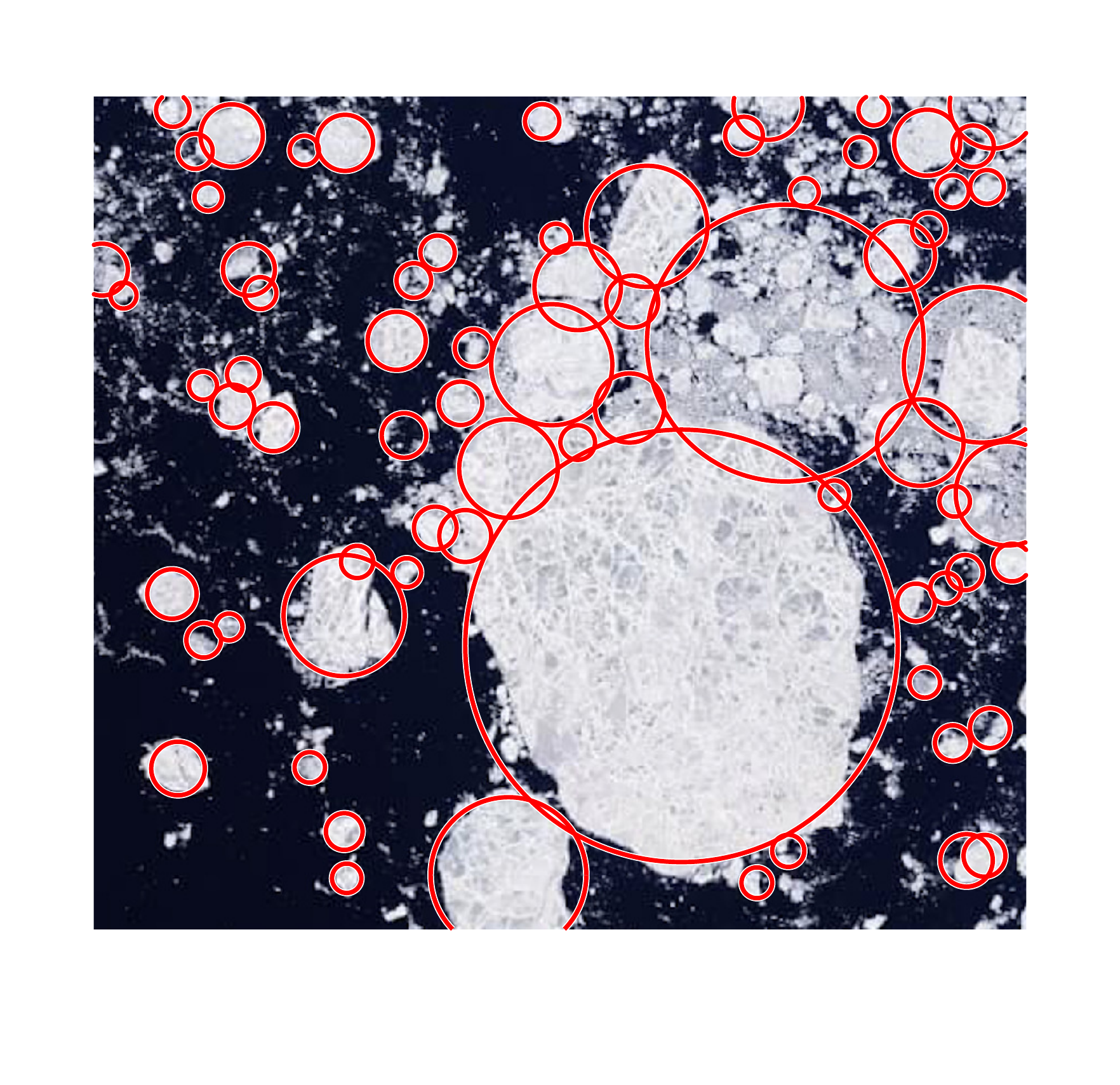}
	\includegraphics[width=2.5in]{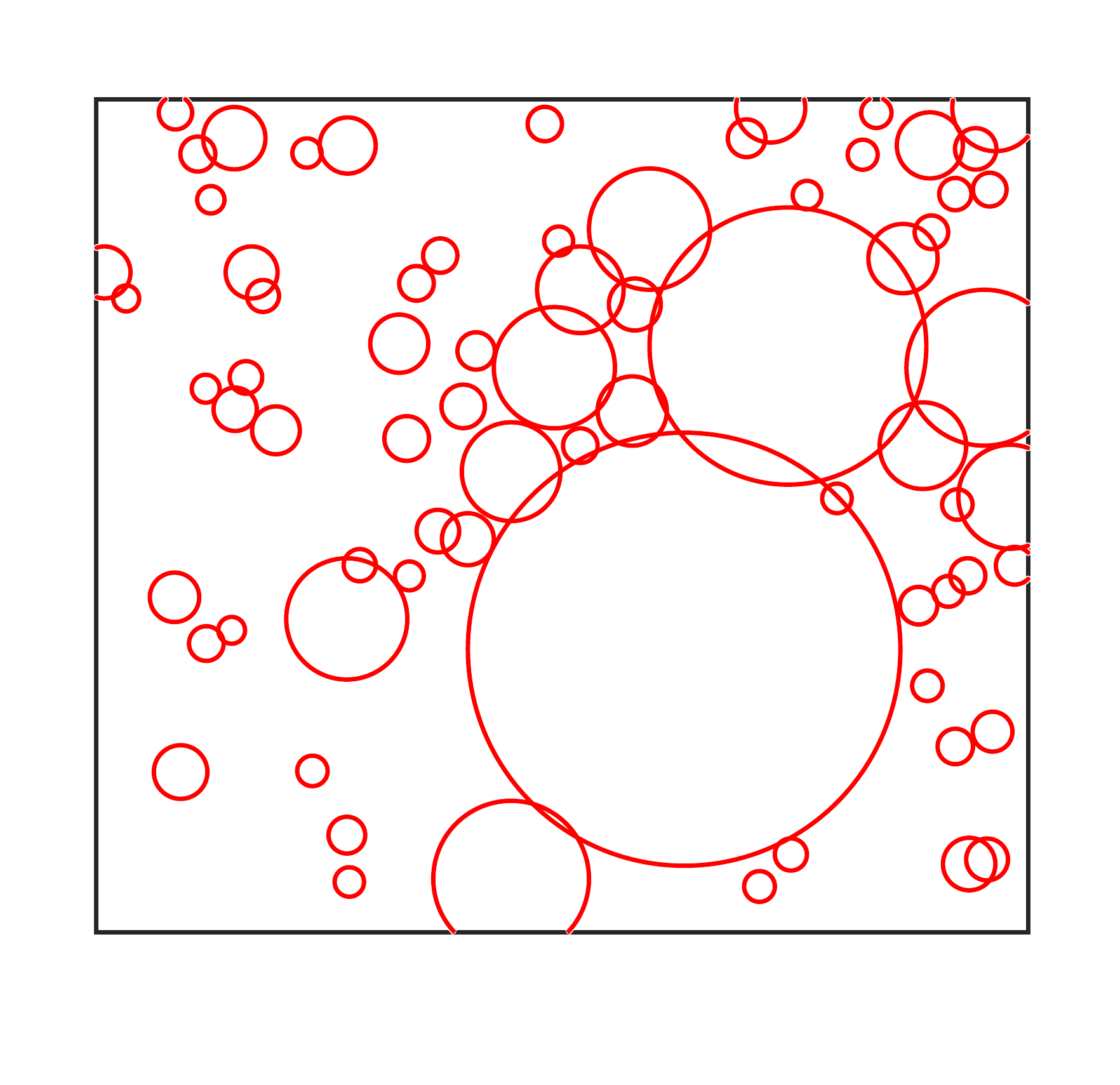}
	\caption{Sea ice floe characterization using cylinders.}
	\label{fig:ofloe}
\end{figure}

Assume there are in total $L$ floes in the system. For simplicity in the mathematical description, all the floes are assumed to be cylinders. Therefore, the geometry of the $l$-th ice floe is determined by the thickness $h^l$ and the horizontal radius $r^l$. The mass is $m^l=\rho_{ice}\pi (r^l)^2 h^l$, where $\rho_{ice}$ is the density of sea ice floes. Assuming the floes all being cylinders is crude but is reasonable, as is illustrated in Figure \ref{fig:ofloe}.

The size and the thickness vary for different floes. Nevertheless, they satisfy certain statistical laws, according to observational data.
The floe size distribution satisfies a power law \cite{stern2018seasonal}
\begin{equation}\label{floe_size_pdf}
p(r) = a\frac{\kappa^a }{r^{a+1}},
\end{equation}
where $r$ is the diameter of the floe, and $\alpha$ and $\kappa$ are parameters. The floe thickness distribution follows a Gamma distribution, whose density function is \cite{bourke1987sea, toppaladoddi2015theory}
\begin{equation}\label{floe_thickness_pdf}
p(h) = \frac{1}{\Gamma(k)\theta^k}h^{k-1}e^{-\frac{h}{\theta}}
\end{equation}
with $k$ and $\theta$ being the shape and scale parameters.. These are both common choices in practice.

\subsection{The equations of motion}
Let superscript $\cdot^l$ denote a state variable of the $l$-th floe with $l=1, \ldots, L$.
The model dynamics contains three sets of state variables:
\begin{enumerate}
  \item [1).] the position $\mathbf{x}^l=(x^l,y^l)^\mathtt{T}$ and the angular location $\Omega^l$ of each floe,
  \item [2).] the velocity $\mathbf{v}^l=(u^l,v^l)^{\mathtt{T}}$ and the angular velocity $\boldsymbol{\omega}^l = \omega^l\hat{\mathbf{z}}$ of each floe, where
 $\hat{\mathbf{z}}$ is the unit vector along the $z$-axis (perpendicular to the $(x,y)$ plane), and
  \item [3).] the ocean surface velocity $\mathbf{u}_o$.
\end{enumerate}

For each floe, Newton's law gives the equations of motion
\begin{equation}\label{eq:xv}
\frac{\d\mathbf{x}^l}{\d t} = \mathbf{v}^l\, ,
\qquad
m^l\frac{\d\mathbf{v}^l}{\d t} = \mathbf{f}^l_o+\mathbf{f}^l_c\, ,
\end{equation}
where the total force involves contributions $\mathbf{f}^l_o$ and $ \mathbf{f}^l_c$ induced by ocean drag forces and floe contact forces, respectively. The ocean drag force obeys the quadratic law \cite{herman2016discrete, damsgaard2018application},
\begin{equation}\label{eq:vo}
  \mathbf{f}^l_o = \widetilde\alpha^l(\mathbf{u}_o -\mathbf{v}^l)|\mathbf{u}_o -\mathbf{v}^l|,
\end{equation}
where $\widetilde\alpha^l= d_o\rho_o \pi  (r^l)^2$ with the constant $d_o$ being the ocean drag coefficient and $\rho_o$ being the density of ocean water. The floe contact force can be further partitioned as
\begin{equation} \label{eq:vc}
    \mathbf{f}_c^l = \sum_{j}(\mathbf{f}_\mathbf{n}^{lj} + \mathbf{f}_\mathbf{t}^{lj}) := \sum_{j}({f}_\mathbf{n}^{lj}\cdot\mathbf{n}^{lj} +  {f}_\mathbf{t}^{lj}\cdot\mathbf{t}^{lj}),
\end{equation}
where the contact force consists of the normal and tangential components, $\mathbf{f}_\mathbf{n}^{lj}$ and $\mathbf{f}_\mathbf{t}^{lj}$ with $\mathbf{n}$ and $\mathbf{t}$ being the unit vectors along the normal and the tangential directions, respectively. The superscript $\cdot^{lj}$ specifies the force from the $j$-the floe to the $l$-th one. The contact force is nonzero when two floes are in contact with each other, i.e., $\delta_n^{lj} \equiv d^{lj} - (r^l+r^j) <0$, where $d^{lj} = |\mathbf{x}^l-\mathbf{x}^j|$ represents the distance between $\mathbf{x}^l$ and $\mathbf{x}^j$. The normal force $\mathbf{f}_\mathbf{n}^{lj}$ is a resistive force to axial compressive stress between two cylindrical ice floes. This force satisfies Hooke's linear elasticity law, i.e., $\mathbf{f}_\mathbf{n}^{lj} = c^{lj}E^{lj}\delta_n^{lj}\mathbf{n}^{lj}$, where $E^{lj}$ is Young's modulus and $c^{lj}$ is the chord length in the transverse direction of the cross-sectional area. The tangential force represents the resistance against slip between floes by limiting relative tangential movement \cite{cundall1979discrete}. That is, $\mathbf{f}_\mathbf{t}^{lj} = c^{lj}G^{lj}v_\mathbf{t}^{lj}\mathbf{t}^{lj}$,
where $G^{lj}$ is the shear modulus and $v_\mathbf{t}^{lj} = \big[(\mathbf{v}^j + \boldsymbol{\omega}^j \times \mathbf{r}^j) - (\mathbf{v}^l + \boldsymbol{\omega}^l \times \mathbf{r}^l)\big]\cdot\mathbf{t}^{lj}$ with $\mathbf{r}^j$ and $\mathbf{r}^l$ being the radius multiplied by the associated normal vector. The normal direction of $\mathbf{r}^j$ is defined by pointing towards the center of the $l$-th floe while $\mathbf{r}^l$ goes the opposite direction. It is important to note that the Coulomb friction law is also used, and it plays an important role in limiting the tangential force relative to the magnitude of the normal force \cite{hopkins2004discrete}. That is, $|\mathbf{f}_\mathbf{t}^{lj}| \leq \mu^{lj}|\mathbf{f}_\mathbf{n}^{lj}|$, where $\mu^{lj}$ is the coefficient of friction that characterizes the condition of the surfaces of the two floes in contact.

On the other hand, the angular velocity $\omega^l$ is given by the rate of change of the angular position $\Omega^l$ in time,
\begin{equation}
  \frac{\d\Omega^l}{\d t} = \omega^l.
\end{equation}
The governing equation of the angular velocity is
\begin{equation}\label{eq:av}
  I^l\frac{\d\omega^l\hat{\mathbf{z}}}{\d t} = \sum_j (r^l\mathbf{n}^{lj}\times \mathbf{f}_\mathbf{t}^{lj}) + t_o^l\hat{\mathbf{z}},
\end{equation}
where $I^l=m^l(r^l)^2$ is the moment of inertia. The first term on the right hand side of \eqref{eq:av} comes from the torque induced by the contact forces while the second part is the torque from ocean drag. The torque is given by
\begin{equation}\label{eq:otk}
  t_o^l\hat{\mathbf{z}} = \widetilde{\beta}^l \left(\frac{\nabla\times \mathbf{u}_o}{2}-\omega^l\hat{\mathbf{z}}\right)\left|\frac{\nabla\times \mathbf{u}_o}{2}-\omega^l\hat{\mathbf{z}}\right|,
\end{equation}
where $\widetilde{\beta}^l = d_o \rho_o\pi (r^l)^4$.

In reality, ocean and sea ice are two-way coupled. Some of the main feedback of sea ice to the ocean comes from the change of the sea surface temperature and the salinity due to the melting and freezing of the sea ice floes. Vice versa, the surface temperature and salinity also affect the melting and freezing of the sea ice floes. They are usually modeled as time-dependent diffusion-advection equations, which bring high computational costs to solve. Since the two-way coupling involves many more complicated modeling procedures and is not directly related to the superfloe parameterization, only a one-way coupling from ocean to sea ice floes is adopted here for simplicity. Such a one-way coupling still allows studying the DA skill as the ocean drag force is one of the dominant contributors to the sea ice motion. Another advantage of adopting the simplified coupling is that it facilitates the use of a general spectrum representation for describing the ocean dynamics.
A set of linear stochastic models is utilized to model each Fourier mode $\hat{u}_{\mathbf{k},\zeta}$ of the ocean current, where the index $\mathbf{k}=(k_1,k_2)$ represents the two-dimensional Fourier wavenumber and the index $\zeta$ is an indicator for different types of the modes (such as the geophysically balanced and the unbalanced ones) associated with the same wavenumber that characterize the ocean flow field. Note that it is well understood that if the underlying ocean model of $\mathbf{u}_o$ is nonlinear, then the time evolution of each Fourier coefficient $\hat{u}_{\mathbf{k},\zeta}$ is driven by a nonlinear deterministic process. Nevertheless, a linear model with additional stochastic noise is a widely used and reasonable representation to approximate the nonlinear deterministic time evolution of $\hat{u}_{\mathbf{k},\zeta}$, especially as a forecast model of DA \cite{majda2016introduction,farrell1993stochastic,berner2017stochastic,branicki2018accuracy,majda2018model,li2018spatial,li2020predictability}. The fundamental mechanism of such an approximation is to stochastically parameterize the effect of the nonlinearity by random noise, which allows an effective quantification of the uncertainty for the underlying ocean dynamics that is required in DA. The governing equation of $\hat{u}_{\mathbf{k},\zeta}$ is given by
\begin{equation}\label{eq:ofm}
  \frac{\d\hat{u}_{\mathbf{k},\zeta}}{\d t} = \big((-d_{\mathbf{k},\zeta} + i\phi_{\mathbf{k},\zeta})\hat{u}_{\mathbf{k},\zeta} + f_{\mathbf{k},\zeta}\big) + \sigma_{\mathbf{k},\zeta}\d \dot{W}_{\mathbf{k},\zeta},
\end{equation}
where $d_{\mathbf{k},\zeta}$, $\phi_{\mathbf{k},\zeta}$ and $\sigma_{\mathbf{k},\zeta}$ are real numbers, representing the damping coefficient, the phase speed of the associated waves and the strength of the stochastic forcing, respectively. On the other hand,
$\dot{W}_{\mathbf{k},\zeta}$ is a complex-valued white noise while $f_{\mathbf{k},\zeta}(t)$ stands for the large-scale deterministic forcing, which is also complex-valued. In practice, the stochastic model is calibrated from the output of the original nonlinear physical model to guarantee the skill of capturing the statistical and basic dynamical behavior. In particular, the parameters of the stochastic model are determined from several key statistics, including the mean, variance and decorrelation time, that are computed from the original nonlinear physical model.  
Define a vector $\mathbf{\hat{u}}_o$ that collects all $\hat{u}_{\mathbf{k},\zeta}$ for different $\mathbf{k}$ and $\zeta$, the spectrum representation of the ocean dynamics can be written into a concise form as
\begin{equation}\label{eq:ofmv}
  \d\mathbf{\hat{u}}_o = \big(\mathbf{L}_\mathbf{u}\mathbf{\hat{u}}_o + \mathbf{F}_\mathbf{u}\big)\d t + \boldsymbol\Sigma_\mathbf{u}\d \mathbf{W}_\mathbf{u}.
\end{equation}
Applying an inverse Fourier transform, the ocean velocity in the physical space can be reconstructed as
\begin{equation}\label{eq:G}
  \mathbf{u}_o = \mathbf{G}(\mathbf{x})\mathbf{\hat{u}}_o,
\end{equation}
where $\mathbf{G}(\mathbf{x})$ is the inverse Fourier transformation matrix.

\subsection{Summary}

Summarizing the above governing equations, the coupled ocean-sea ice system is
\begin{subequations}\label{eq:pde}
\begin{align}
  \frac{\d\mathbf{x}^l}{\d t} &= \mathbf{v}^l,\label{pde_x}\\
  \frac{\d\Omega^l}{\d t}  &= \omega^l,\label{pde_Omega}\\
  m^l\frac{\d\mathbf{v}^l}{\d t} &=    \sum_{j=1}^L(\mathbf{f}_\mathbf{n}^{lj} + \mathbf{f}_\mathbf{t}^{lj}) + \widetilde{\alpha}^l\left(\mathbf{G}(\mathbf{x}^l) \hat{\mathbf{u}}_o -\mathbf{v}^l\right)\left|\mathbf{G}(\mathbf{x}^l) \hat{\mathbf{u}}_o -\mathbf{v}^l\right|,\label{pde_v}\\
  I^l\frac{\d\omega^l}{\d t} &= \sum_{j=1}^L ( r^l\mathbf{n}^{lj}\times \mathbf{f}_\mathbf{t}^{lj}) \cdot\hat{\mathbf{z}} + \widetilde{\beta}^l \left(\nabla\times \mathbf{u}_o/2-\omega^l\hat{\mathbf{z}}\right)\left|\nabla\times \mathbf{u}_o/2-\omega^l\hat{\mathbf{z}}\right|,\label{pde_omega}\\
  \frac{\d\mathbf{\hat{u}}_o}{\d t} &= \big(\mathbf{L}_\mathbf{u}\mathbf{\hat{u}}_o + \mathbf{F}_\mathbf{u}\big) + \boldsymbol\Sigma_\mathbf{u} \dot{\mathbf{W}}_\mathbf{u}(t)\label{pde_u},
\end{align}
\end{subequations}
where $l=1,2, \cdots, L$. 

Despite being simplified compared with the operational models of the sea ice, the coupled system captures many key features of the sea ice floe dynamics. Note that  \eqref{eq:pde} is highly nonlinear due to the quadratic terms in the linear and angular momentum equations as well as the nonlinear operator $\mathbf{G}(\mathbf{x}^l)$, which is an exponential function of $\mathbf{x}^l$. The coupled system is also a high-dimensional system. The total dimension of the system is $6L+D^o$, where $D^o$ is the number of degrees of freedom of the ocean.


\section{An efficient physics-constrained superfloe parameterization} \label{sec:sp}

One of the most computationally challenging aspects of DEM sea ice simulations is the high dimensionality of the system. In fact, a typical operational DEM model consists of at least a few tens of thousands of floes. Therefore, effective parameterizations for these Lagrangian DEM models are crucial for improving computational efficiency while retaining the key dynamical features. Suitable parameterizations also facilitate efficient DA. It is important to note that the reduced-order system by simply removing the small-scale floes, as an analog to the bare truncation in typical turbulent systems, often brings about a large error since the interactions from small-scale floes to the large-scale ones via contact forces have a significant contribution to the overall dynamics \cite{tuhkuri2018review}.

\subsection{Key features of the superfloe parameterization}
\begin{figure}[h!]
\centering
	\includegraphics[width=5in]{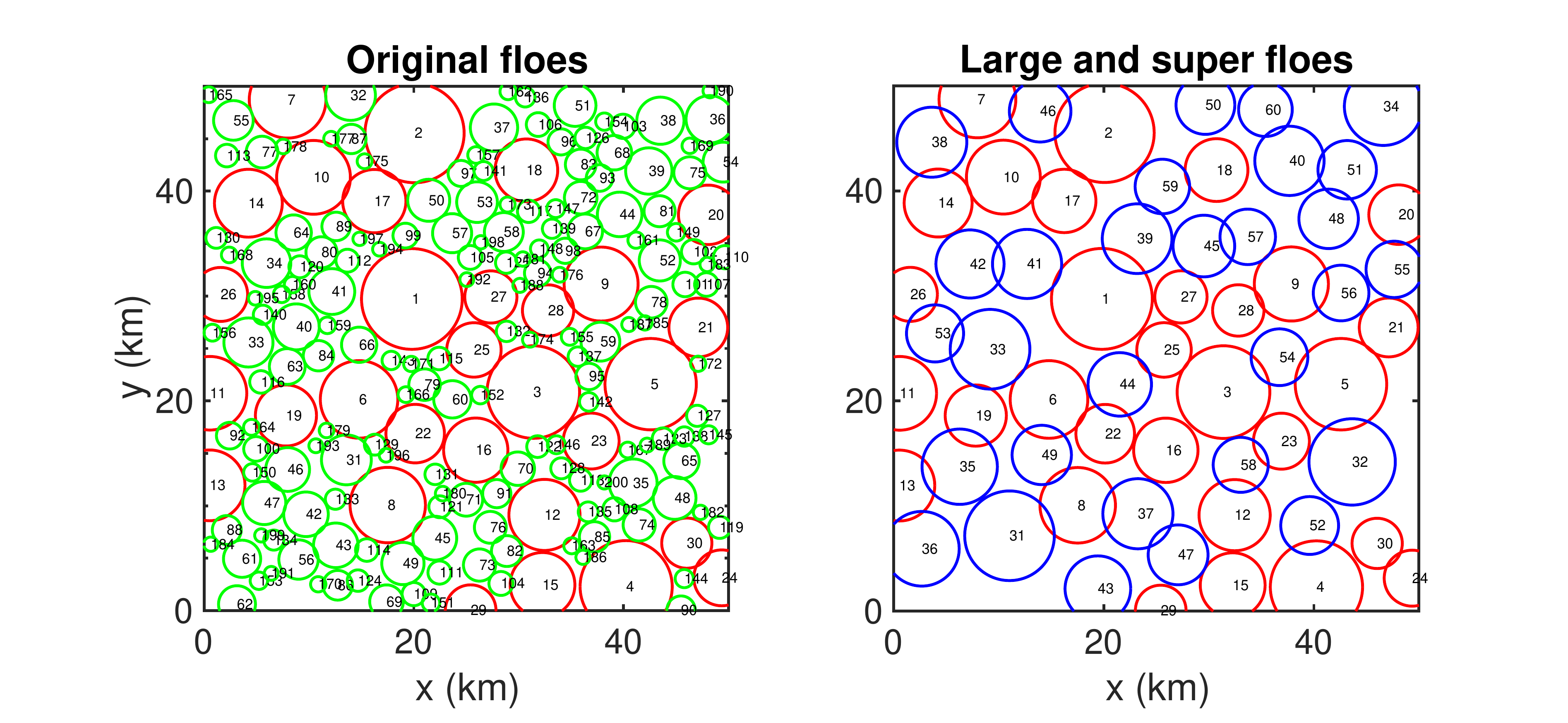}
	\caption{An example of superfloe parameterization. There are in total $200$ floes in the original system.
	Red: $L_0=30$ largest floes; Green: $L-L_0=170$ small floes; Blue: $L_s=30$ superfloes.}
	\label{fig:200to60}
\end{figure}

The new parameterization developed here exploits artificial sea ice floes, which are named as ``superfloes'', to act as a substitute for the small-scale floes in the reduced-order system. Each superfloe aims at approximating a cluster of the small-scale floes. Therefore, only a small number of superfloes is sufficient to effectively parameterize all the small-scale floes, which significantly reduces the dimension of the resulting system. See Figure~\ref{fig:200to60} for an illustration. Different from many parameterizations that involve empirical tunings, physics constraints are systematically incorporated into the new superfloe parameterization. Specifically, the superfloe parameterized system is designed to preserve the following physical quantities in the original system. These quantities are:
\begin{enumerate}
  \item [1).] the mass,
  \item [2).] the concentration,
  \item [3).] the linear momentum, and
  \item [4).] the angular momentum.
\end{enumerate}
Including the mass constraint is natural, which is also the basis for retaining many other quantities, such as the concentration. On the other hand, the momentum is a more robust quantity to utilize than the energy as a physics constraint. In fact, the energy transfer between floes is very complicated especially in the presence of collisions, where part of the energy is dissipated. In contrast, the total momentum is conserved instantaneously when the collision occurs and the momentum is simply interchanged between different floes.

Assume that there are $L$ floes, which are sorted in ascending order according to their floe size (i.e., the radius here). The first $L_0$ large-scale floes are retained in the parameterized system, where $L_0\ll L$, while the remaining $L-L_0$ relatively small floes are parameterized by $L_s$ superfloes, where $L_s \ll L- L_0$. Therefore, there are only in total $L_r = L_0 + L_s$ floes in the parameterized system. Then the constraints of the mass, the linear momentum, and the angular momentum are given by

\begin{equation} \label{eq:menew}
\begin{aligned}
\mbox{Mass:}&&\quad {m}_{total} =&  \sum_{j=1}^{L_0}  m^j + \sum_{j=L_0+1}^{L} m^j &=& \sum_{j=1}^{L_0}  m^j + \sum_{k=L_0+1}^{L_r} \tilde{m}^k, \\
\mbox{Concentration:}&&\quad {c}_{total} =&  \sum_{j=1}^{L_0}  (r^j)^2 + \sum_{j=L_0+1}^{L} (r^j)^2 &=& \sum_{j=1}^{L_0}  (r^j)^2 + \sum_{k=L_0+1}^{L_r} (\tilde{r}^k)^2, \\
\mbox{Linear momentum:}&&\quad {\mathbf{p}}_{total} =&   \sum_{j=1}^{L_0}  m^j\mathbf{v}^j + \sum_{j=L_0+1}^{L} m^j\mathbf{v}^j &=& \sum_{j=1}^{L_0}  m^j \mathbf{v}^j+ \sum_{k=L_0+1}^{L_r} \tilde{m}^k\tilde{\mathbf{v}}^k, \\
\mbox{Angular momentum:}&&\quad {\mathbf{L}}_{total} =&  \sum_{j=1}^{L_0}  I^j  {\mathbf{\omega}}^j + \sum_{j=L_0+1}^{L} I^j  {\mathbf{\omega}}^j &=& \sum_{j=1}^{L_0}  I^j  {\mathbf{\omega}}^j +  \sum_{k=L_0+1}^{L_r}  \tilde I^k \tilde{\mathbf{\omega}}^k,
\end{aligned}
\end{equation}
where $\tilde{\cdot}$ denotes the quantities associated with the superfloe parameterization, and the constants in the expression of the concentration have been ignored.

\subsection{A superfloe parameterization algorithm} 
The $k$-th superfloe  combines a cluster of $J$ nearby small-scale floes into one superfloe. The mass constraint of the superfloe leads to
\begin{equation} \label{eq:spfm}
\tilde{m}^{k} = \sum_{j=1}^J m^j.
\end{equation}
Similarly, the area of the superfloe equals the sum of the areas of the $J$ small-scale floes, which guarantees the constraint of the sea ice concentration. For a cylinder floe, the area is determined by its radius, i.e., $\pi r^2$.
Thus, the radius of the $k$-th superfloe  is
\begin{equation} \label{eq:spfr}
\tilde{r}^{k} = \sqrt{\sum_{j=1}^J (r^j)^2}.
\end{equation}
Now with the expressions of the mass and the area in hand, the thickness of the $k$-th superfloe can be calculated
\begin{equation} \label{eq:spfh}
\tilde{h}^{k} = \frac{\tilde{m}^{k}}{ \rho_{\text{ice}} \pi^2 (\tilde{r}^{k})^2 }.
\end{equation}
On the other hand, the position of the superfloe is given by the center of the mass of the $J$ small-scale floes,
\begin{equation} \label{eq:spfx}
\tilde{\mathbf{x}}^{k} = \frac{1}{\tilde{m}^k} \sum_{j=1}^J m^j \mathbf{x}^j.
\end{equation}
Next, the velocity of the $k$-th superfloe  is calculated from the constraint of the linear momentum,
\begin{equation} \label{eq:spfv}
\tilde{\mathbf{v}}^{k} = \frac{1}{\tilde{m}^{k}}\sum_{j=1}^J m^j \mathbf{v}^j.
\end{equation}
Similarly,  the angular velocity of the $k$-th superfloe is given by the constraint of the angular momentum as
\begin{equation} \label{eq:spfw}
\tilde{\omega} = \frac{1}{\tilde{I}^{k}}\sum_{j=1}^J I^j \omega^j,
\end{equation}
where the moment of inertia of the superfloe is $\tilde{I}^{k} = \tilde{m}^{k}  (\tilde{r}^{k})^2$.

The superfloe parameterization is summarized in Algorithm \ref{alg:fnr}.
\begin{algorithm}[h!] 
\rule{14.2cm}{0.05cm}
\begin{algorithmic}
\STATE{Initialize the system of $L$ total floes and set floe number target $L_0, L_s, L_r=L_0+L_s$.}
\STATE{Keep the $L_0$ largest floes.}
\WHILE{$L_r>L_0+L_s$}
\STATE{Sort all the small floes and superfloes in descending order with respect to their radii.}
\STATE{Start from the smallest floe to seek a group of neighboring floes.}
\IF{the distance between the smallest floe and nearby floes is large}
\STATE{Delete the smallest floe as it is well-isolated.}
\STATE{Update the number of floes $L_r$.}
\ELSE
\STATE{Create a superfloe with quantities evaluated by equations \eqref{eq:spfm}--\eqref{eq:spfw}.}
\STATE{Update the number of floes $L_r$.}
\ENDIF
\ENDWHILE
\STATE{Return a new set of floes.}
\end{algorithmic}
\rule{14.2cm}{0.05cm} \\[0.3cm]
\caption{Superfloe parameterization}
\label{alg:fnr}
\end{algorithm}

In general, a small floe has a small contact force and hence a small impact on its neighbouring floes. If a small floe is far away from all other nearby floes, it has no contact force to its nearby floes, at least for a short period of time when it is not interacting with other floes. We hence delete those floes as in Algorithm 3.1.
In our numerical experiments, we set this distance threshold as $\sqrt{2}(r_{\text{smallest}} + r_{\text{nearby}})$ where $r_{\text{smallest}}$ is the current smallest floe and $r_{\text{nearby}}$ is the floe closest to this smallest floe.
This leads to a parameterized system that satisfies the mass conservation law in an approximate way. Alternatively, one can keep it in the system to retain exact mass conservation. If keeping it, then we slightly modify the algorithm and search the nearby floes of the next smallest floes. In real applications, such isolated and small floes are very rare due to the collision, melting, and welding features of sea ice. Thus, we simply remove them as in Algorithm 3.1 for simplicity.
Lastly, to identify groups of neighboring floes,
while other algorithms are possible,
in our numerical experiments, we propose to group gradually two nearby floes where one of the floes is the smallest floe in the current floe setting state.
This technique has the advantage of a lower computational cost in the sense that there is no need to determine the size (floe number) of a group of neighboring floes and to determine the group a floe belongs to.
Another advantage is that it avoids a superfloe becoming too large compared to other floes as a DEM system with larger floe size ranges requires smaller time-marching step sizes.

For the dynamical equations of motion, the superfloe parameterization
has another advantageous feature: the same dynamical equations from \eqref{eq:pde}
can be used for superfloes and for ordinary floes.
As a result, no additional specifications are needed for
the interactions between superfloes and ordinary floes,
nor for the interactions between one superfloe and another superfloe,
since they all interact in the same way that ordinary floes
interact with each other. As one possible modification, one might
suspect that a superfloe should have a reduced value of the
Young's modulus in comparison to an ordinary floe, since
a superfloe is less like solid ice than an ordinary floe,
and a superfloe should perhaps have a weaker response upon
a collision. It would be interesting to consider such possibilities
in the future. Here, for simplicity, we investigate the
use of the same Young's modulus in what follows, and we find
that it yields satisfying results.

 %

\subsection{Example of superfloes and computational savings}
Figure \ref{fig:200to60} compares one snapshot of the full floe field and the one with superfloe parameterization. The full floe field contains $200$ floes, while the one with superfloe parameterization retains the largest $30$ floes and parameterizes the other $170$ floes by $30$ superfloes.
It is clear that the groups of the neighbouring floes are reasonably well represented by the superfloes. For instance, the small floes with numbers 55, 77, 113, 165, and 178 are neighbouring small-floes, which are parameterized as a superfloe with number 38 (see the top-left corner of plots in Figure \ref{fig:200to60}).
We also observe some small overlapping of the floes in the right panel of the Figure \ref{fig:200to60}. This is mainly due to the initialization that the floes are dense and the assumption that the floes are cylinders. The overlapped areas can be reduced with a different initialization or floe geometry. Nevertheless, it is worthwhile to point out that such overlaps are also part of the modeling process in many more realistic DEM models. The overlaps are utilized to effectively characterize the  elastic deformation due to  tension or compression of the sea ice floes, which are often assumed to be a rigid body.  

Table \ref{tab:fnr} compares the floe statistics in the full system and those in the reduced-order system with superfloe parameterization. Different rows show the cases with different numbers of the floes $L$ in the full system. It is clear that the reduced-order system with the superfloe parameterization results in the same concentration and the minimum and maximum of the thickness as those in the full system. The minimum radius $r_{min}$ in the superfloe parameterized system is larger than that in the full system since the small-scale floes have been substituted by the superfloes. The maximal radius $r_{max}$ can be larger when a large floe is grouped with its neighbouring floes. These results also indicate the robustness of the superfloe parameterization.

\begin{table}[ht]
\centering
\begin{tabular}{| c cc | ccccc | ccccc | cc |}
\hline
\multicolumn{3}{|c|}{ }  & \multicolumn{5}{c|}{Full system} & \multicolumn{5}{c|}{Superfloe parameterized system} \\[0.1cm] \hline
$L$ & $L_s$ & $L_r$ & $c$ & $r_{\min}$ &  $r_{\max}$ & $h_{\min}$  & $h_{\max}$  & $c$ & $r_{\min}$ &  $r_{\max}$ & $ h_{\min} $  & $ h_{\max} $  \\[0.1cm] \hline
40 & 20 & 40 &   0.34 &   1.51 & 4.22 & 0.19 & 2.54 & 0.34 &   1.51 & 4.22 & 0.19 & 2.54 \\[0.1cm]
60 & 20 & 40 &   0.50 &   1.51&    4.38 &    0.18 &    2.54 & 0.48 &    2.15 &4.38&    0.20 &    2.54 \\[0.1cm]
80 & 20 & 40 &   0.58    & 1.51&    3.98    &0.18&    2.54 & 0.56&   2.71 &   4.06&    0.20&   2.54  \\[0.1cm]
100 & 20 & 40 &   0.75 &    1.51&    4.06 &   0.17& 2.54 & 0.74&    2.34 &    5.25&    0.20&   2.54 \\[0.1cm]
100 & 30 & 60 &   0.75&    1.51&    4.06&   0.17&    2.54&   0.75&    2.31 &    4.46&    0.19&    2.54 \\[0.1cm]
200 & 30 & 60 &   0.78&   0.80&    3.82&    0.17&    3.33& 0.78&    2.47 &    4.46&    0.19&    2.54 \\[0.1cm]\hline
 \end{tabular}
\caption{Comparisons of the floe statistics in the full system and those in the reduced-order system with superfloe parameterization. The statistics include the sea ice concentration $c$, the minimum and maximum of the radius $r_{\min}$ and $r_{\max}$, and the minimum and maximum of the thickness $h_{\min}$ and $h_{\max}$. Recall that $L$ is the total number of the floes in the full system,
$L_s$ is the number of superfloes, and
$L_r$ is the total number of the floes in the reduced-order system with the superfloe parameterization.
Radius unit: km; thickness unit: m.}
\label{tab:fnr}
\end{table}

\begin{figure}[h!]
	\centering
	\includegraphics[width=4in]{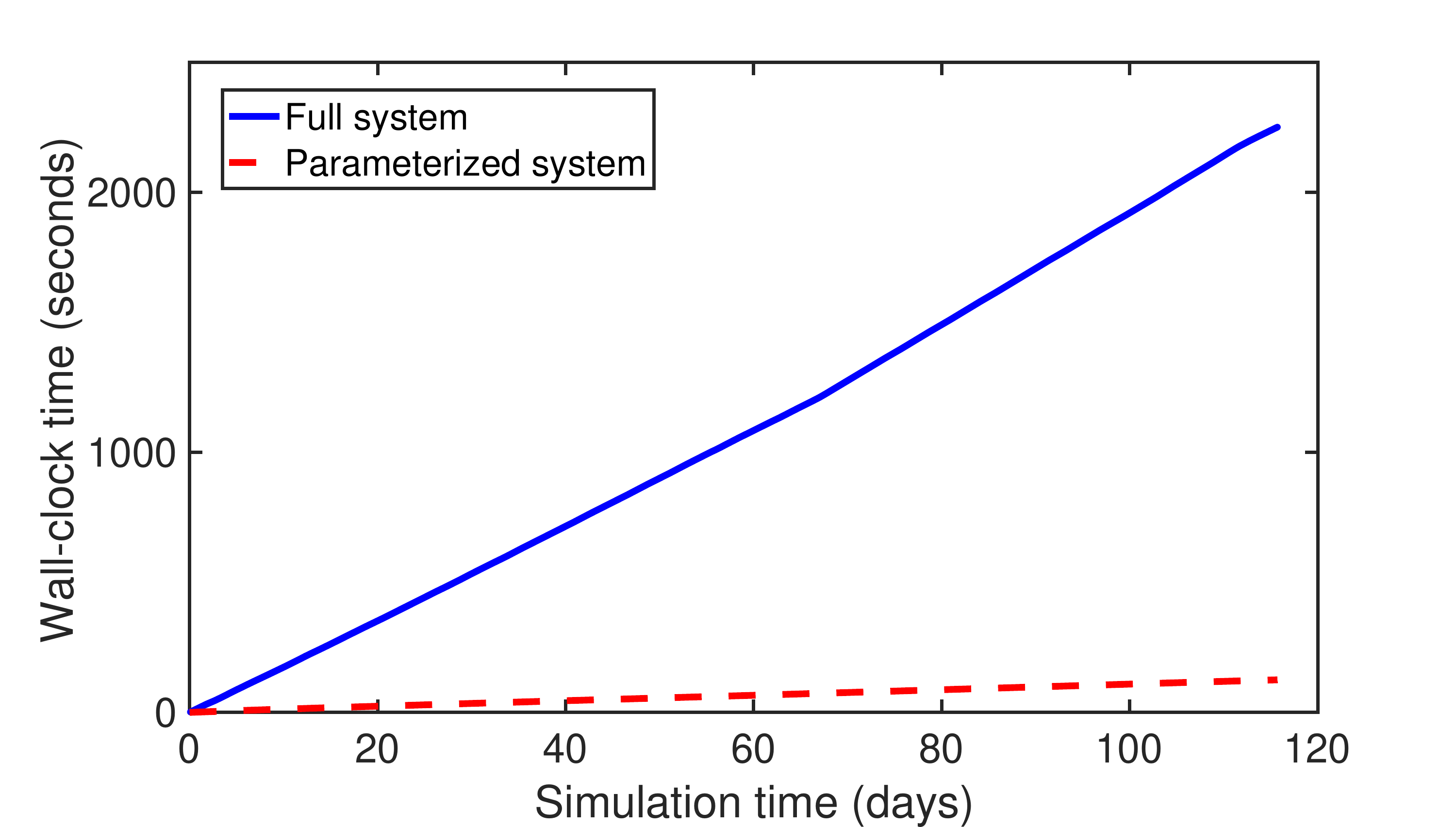}
	\caption{Computational time cost comparison of the original (200 floes) and superfloe parameterized (30 large floes and 30 superfloes) systems.}
	\label{fig:time}
\end{figure}

As a brief, first look at evolutionary simulations, Figure \ref{fig:time} compares the computational cost of the evolution of the full system \eqref{eq:pde} and that of the reduced-order system with the superfloe parameterization. The final time of the simulation is about $T=120$ days.  The parameterization significantly reduces the simulation time costs. Moreover, since there are fewer floes in the parameterized system, the computational storage costs are also reduced.


\section{Uncertainty quantification (UQ) with superfloe parameterization} \label{sec:uq}
Now we consider, in more detail, the superfloe parameterization and its incorporation into the floe model \eqref{eq:pde} to form a reduced-order system, the simulation of which will be compared with the full system \eqref{eq:pde}. The primary goal here is to understand the uncertainty quantification skill of the reduced system with the superfloe parameterization. The uncertainty quantification is two-fold. On one hand, we expect the reduced system with the superfloe parameterization captures the time evolution of the uncertainty as the full system for a short-term forecast. On the other hand, it is important to see if the reduced system is able to recover the long-term statistics especially the non-Gaussian distribution with fat tails. Throughout this paper, the parameters in the floe size distribution \eqref{floe_size_pdf} are $\alpha=1$ and $\kappa=1.5$ while those in the thickness distribution \eqref{floe_thickness_pdf} are $k=2$ and $\theta=1.3$. A square domain of scale $50$km$\times50$km is used here, mimicking the marginal ice zone. The ocean field is generated from a truncated linear shallow water system \cite{majda2003introduction, vallis2017atmospheric}. It contains 26 Fourier modes, with 8 geophysically balanced (GB) modes and 18 gravity modes. The GB modes are incompressible and they are slowly varying in time. On the other hand, the gravity modes are compressible and  have fast oscillations. The Rossby number is Ro$=0.1$ such that the gravity modes lie in a much faster time scale than the GB modes. The damping coefficients for all the Fourier modes in \eqref{pde_u} are $0.5$. The noise coefficients of the GB modes are all $0.1$ while those of the gravity modes are all $0.05$. There is no deterministic forcing in the ocean equation. These parameters allow the energy in the GB part of the flow to be roughly twice as much as that in the unbalanced gravity modes. The ocean velocity is of order $0.1$m/s, which is consistent with observations. The numerical integration time step is $\Delta{t} = 25$ seconds to resolve the gravity modes. The other parameters as well as their physical units are listed in Table \ref{tab:param}.

\subsection{Short-term behavior of the reduced-order system with the superfloe parameterization}

\begin{figure}[ht]
\centering
	\includegraphics[width=6in]{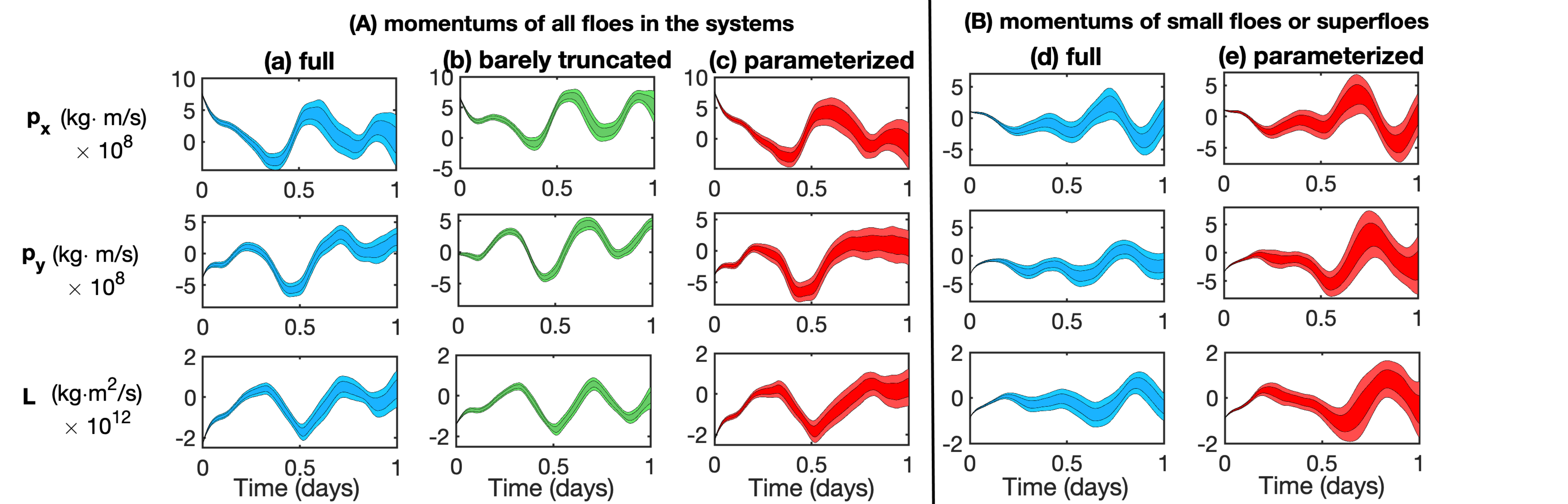}
	\caption{Comparison of the short-term behavior of the momentums.  All three systems start with the same initial condition. The uncertainty increases as the systems run forward in time due to the random forcing in the systems. The dark and light shading areas show the ensemble spread corresponding to 1 and 2 standard deviations of the ensembles, respectively.}
	\label{fig:mom}
\end{figure}

We start with studying the short-term dynamics of the superpfloe parameterized system.
Figure \ref{fig:mom} shows the time evolution of the momentum. In the experiment here, there are in total $18$ floes in the full system. Only the $6$ largest floes are retained in the bare truncation system. On the other hand, the superfloe parameterized system contains the $6$ largest floes and $6$ superfloes.
The momentum variables are collected component-wise as sums over all floes in each system.  All three systems start with the same initial condition. An ensemble forecast with 1000 ensemble members is used for the short-term forecast here. The initial values of the ensembles are all assigned to be the same, which equals  the exact value from the full model. The uncertainty increases as the systems run forward in time due to the random forcing in the systems. The goal here is to compare the time evolution of the uncertainty in the ensemble forecast using the three systems.

Panel (A) shows the total momentum of all the floes in each system. It is clear that the time evolution of the momentum, as well as the associated uncertainty in the original system, are well captured by the use of superfloes due to the physics constraints in the superfloe parameterization. In contrast, the uncertainty in the bare truncation model is severely underestimated, which indicates the necessity of parameterizing the effects from the small-scale floes using the superfloe parameterization. In fact, Panel (B) illustrates the total momentum after removing the 6 large floes in the original and parameterized system, where Panel (d) shows the total momentums of the $12$ small floes in the original system while Panel (e) shows the total momentum of the $6$ superfloes in the parameterized system. Such a comparison implies that the superfloes indeed recover the uncertainty propagation in the small-scale floes.


\subsection{Long-term statistics of the reduced-order system with the superfloe parameterization}
The focus of this subsection is on comparing the long-term statistical behavior of the reduced-order system with the superfloe parameterization versus the full system.\medskip

\begin{figure}[h!]
\centering
	\includegraphics[width=5in]{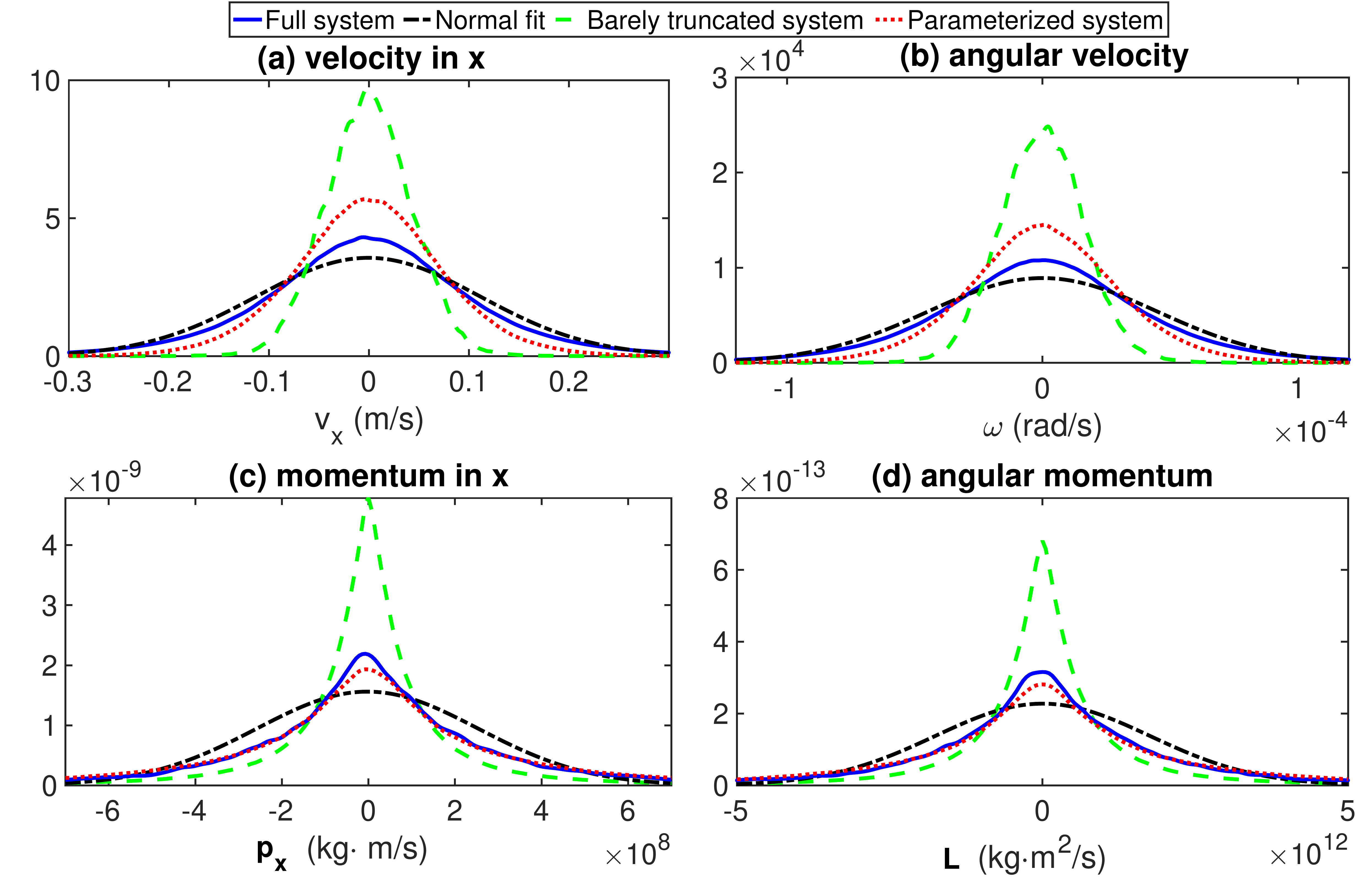}
	\caption{Comparison of the probability density functions (PDFs) of various floe physical quantities.
	There are 30 large floes and 30 superfloes in the parameterized system. }
	\label{fig:uqpdfs}
\end{figure}

\noindent\textbf{(a). Statistics of several key physical quantities.}\\
Figure \ref{fig:uqpdfs} compares the probability density functions (PDFs) of the velocities, angular velocities, linear momentum, and angular momentum associated with the large-scale floes in three different systems:
\begin{enumerate}
  \item [1).] the full system, which contains in total $L=200$ floes;
  \item [2).] the superfloe parameterized system,  where the $L_0=30$ large-scale floes from the full system are retained and the remaining $170$ floes are parameterized by $L_s=30$ superfloes; and
  \item [3).] the bare truncation system, where  only the $L_0=30$ large-scale floes are retained while the other $170$ floes are completely ignored.
\end{enumerate}
The PDFs are based on the simulations over the time interval from $T=30$ (days) to $T=120$ (days).

It is clear that the statistics associated with the superfloe parameterized system resemble those of the full system. In contrast, the barely truncated system has completely different statistical behavior. The main difference between these two approximations is that the superfloes mimic the small-scale floes to provide statistically accurate contact forces to the large-scale floes. These contact forces are important to recover the statistics of all the quantities. Note that despite the PDFs of the velocity and the angular velocity being nearly Gaussian, the momentum and angular momentum have highly non-Gaussian statistics. The superfloe parameterized system succeeds in recovering these fat-tailed PDFs.

\noindent\textbf{(b). Statistics of the contact forces.}\\
Now, we take a detailed look at the skill of the superfloe approximations in recovering the contract forces. We define the contact force from all the small-scale floes or the superfloes to the $k$-th large floe in \eqref{pde_v} and \eqref{pde_omega} as
\begin{equation} \label{eq:cf}
\mathbf{f}_c^k = \sum_{j=L_0+1}^{L_f} (\mathbf{f}_\mathbf{n}^{kj} + \mathbf{f}_\mathbf{t}^{kj}) \quad \text{and} \quad
f_\omega^k = \sum_{j=L_0+1}^{L_f} ( r^k\mathbf{n}^{kj}\times \mathbf{f}_\mathbf{t}^{kj}) \cdot\hat{\mathbf{z}},
\end{equation}
where $L_f = L$ for the original system and $L_f=L_r$ for the parameterized system, and $\mathbf{f}_c^k$ has two components along $x$ or $y$ directions, respectively.
%
The statistics are computed based on a long simulation time.

\begin{figure}[h!]
\centering
	\includegraphics[width=6in]{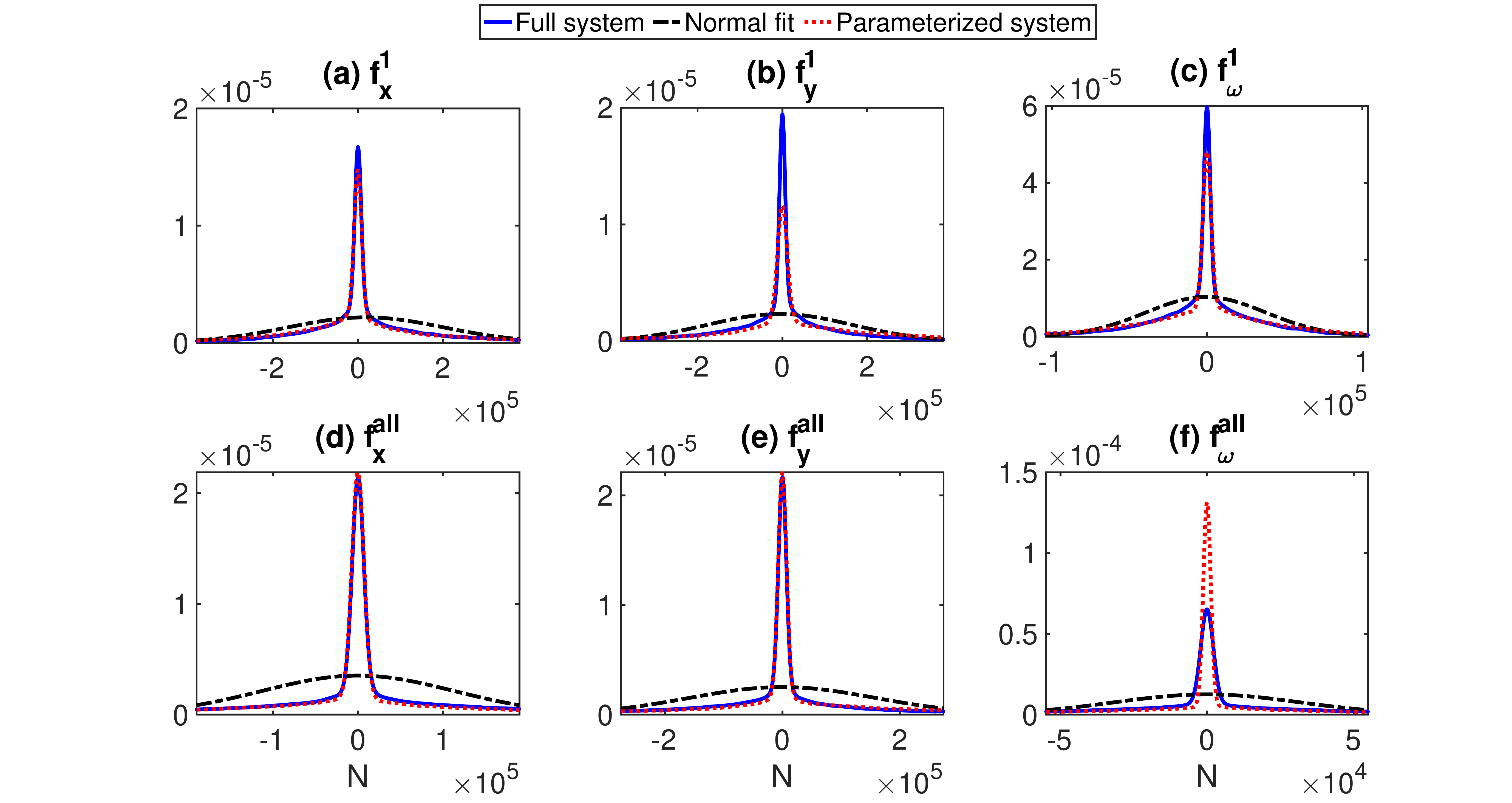}
	\caption{Comparisons of the PDFs  of floe contact forces. There are $18$ floes in total. The $12$ smaller floes are parameterized as $6$ superfloes. Top row: statistics of the largest floe; bottom row: statistics of all the $6$ floes. The normal fit is a fit for the contact forces of the original system.}
	\label{fig:cf}
\end{figure}

The first row of Figure \ref{fig:cf} shows the PDFs of the contact forces to the largest floe $\#1$ while the second row shows those to all the $30$ large-scale floes. Despite a significant dimension reduction of the system, the superfloes succeeds in recovering the highly non-Gaussian statistics of the contact forces.
The results here indicate that the superfloe parameterized system is statistically accurate for describing the features of the large-scale floes with a much reduced computational cost.

\section{Lagrangian data assimilation (DA) with superfloe-based inflation} \label{sec:da}
In addition to approximating the statistical behavior of the large-scale floes in the full system, the superfloe parameterization can also be used to facilitate the Lagrangian DA that recovers the unobserved ocean field by observing the floe trajectories.

In practice, only the large-scale floes are easily identified from the satellite images, which are the observations in the Lagrangian DA. This means the contact force in the equation \eqref{eq:pde} due to the small-scale floes cannot be fully resolved in the forecast model. Simply ignoring the contributions from the small-scale floes (e.g., using the bare truncation system) is expected to have large biases. Noise inflation \cite{anderson2009spatially, li2009simultaneous} is a typical technique that improves the DA skill in the presence of such a model error. However, noise inflation often relies on many ad hoc tunings, which makes it very difficult to apply in practice.

In the following, a systematic noise inflation strategy based on the superfloe parameterization is developed to determine the noise inflation coefficients that significantly advance the DA skill. The ensemble adjustment Kalman filter (EAKF) \cite{anderson2001ensemble} will be utilized as the DA algorithm throughout this section.

\subsection{A superfloe-based noise inflation algorithm} \label{sec:dasp}
Let $\Delta t^{obs}$ be the observational time step and $\Delta t$ be the numerical integration time step. Define $M=\lfloor \Delta t^{obs}/\Delta t \rfloor$ where $\lfloor \cdot \rfloor$ is a floor function.
The superfloe-based noise inflation algorithm is given as follows.

\begin{algorithm}[h!]
\rule{14.5cm}{0.05cm} 
\begin{algorithmic}
\STATE{Develop a superfloe model using Algorithm \ref{alg:fnr} with $L_0$ large floes and $L_s$ superfloes.}
\STATE{Run the superfloe model up to $T=N\Delta t$ and let  $t_j = j \Delta t , j = 0,1,\cdots, N,$.}
\STATE{Store the contact forces of the $l$-th large floe that are from collisions with the $L_s$ superfloes, i.e., $\mathbf{f}_c^l(t_j)$ and $f_\omega^l(t_j)$ in \eqref{eq:cf} for $l=1,2,\cdots,L_0$. }
\FOR{$l=1,2,\cdots,L_0$}
\STATE{Set $f_j = f_\omega^l(t_j), j =1,2,\cdots,N$. }
\STATE{Form a vector $\bfs{F} = (f_1, f_2,\cdots, f_N)^T$. }
\STATE{Set $\bfs{F}_2 = (f_{M+1}, f_{M+2},\cdots, f_N)^T$ and $\bfs{F}_1 = (f_{1}, f_{2},\cdots, f_{N-M})^T$. }
\STATE{Calculate the standard deviation of $\bfs{F}_2-\bfs{F}_1$ and store it as $\tilde{\sigma}^l_\omega$.}
\STATE{Do the same for $\mathbf{f}_c^l(t_j)$ to obtain standard deviation $\tilde{\bfs{\sigma}}^l_c$. }
\ENDFOR
\STATE{Use $\tilde{\bfs{\sigma}}^l_c$ as the noise inflation coefficient in \eqref{pde_v} for each large floe $l, l=1,2,\cdots,L_0$. }
\STATE{Use $\tilde{\sigma}^l_\omega$ as the noise inflation coefficient in \eqref{pde_omega} for each large floe $l, l=1,2,\cdots,L_0$.}
\end{algorithmic}
\rule{14.5cm}{0.05cm} \\[0.3cm]
\caption{Superfloe-based noise inflation}
\label{alg:da}
\end{algorithm}

In Algorithm \ref{alg:da}, the noise inflation is applied to the forward modeling equation by adding $\tilde{\bfs{\sigma}}^l_c \d\mathbf{W}_\bfs{v}^l(t)$ to the right-hand side of \eqref{pde_v} and $\tilde{\sigma}^l_{\omega} \d W_{\omega}^l(t)$ to the right-hand side of \eqref{pde_omega}, respectively.
In the DA setting, we assume that we observe the trajectories of the large floes (partial observation)
and our goal is to estimate the velocities and the Fourier modes of the underlying ocean currents.
To set up the EAKF for the sea ice floe simulation,
we let the state variable $s$ be a vector containing the displacements (locations)
and velocities of all the floes as well as the ocean current Fourier modes.
With this in mind,  the state variable is
\begin{equation*}
s = (x^1_1, \cdots, x^{L_t}_1, x^1_2, \cdots, x^{L_t}_2, \Omega^1, \cdots, \Omega^{L_t}, v^1_1, \cdots, v^{L_t}_1, v^1_2, \cdots, v^{L_t}_2, \omega^1, \cdots, \omega^{L_t}, \mathbf{\hat{u}}_o)^T,
\end{equation*}
where $L_t$ is the total number of floes in the (original or superfloe-parameterized) system.
The Algorithm \ref{alg:da} is then implemented as the core of the forecast model for the EAKF which is based on \cite{anderson2001ensemble}. The observational noise is set to be $\sigma_\textbf{x}=$80m, $\sigma_\Omega = $0.01rad. A summary of the parameters is listed in Table \ref{tab:param}.

The idea here is to compute the averaged strength of the contact force variability over one forecast cycle $\Delta t^{obs}$, utilizing it as the noise inflation coefficient.
Admittedly, the noise inflation coefficient can be computed based on the full model \eqref{eq:pde}. However, while running the idealized model \eqref{eq:pde} is computationally affordable, running a full operational DEM model for a long time is not practical. Since it has been shown in Figure \ref{fig:cf} that the superfloe model succeeds in capturing the highly non-Gaussian statistics of the contact force, the much cheaper superfloe model is more appropriate for determining the noise inflation coefficients in the DA forecast model.

\subsection{Numerical experiments}
The root-mean-square error (RMSE) and the pattern correlation coefficient (PCC) will be utilized as the skill scores to quantify the performance of DA. Denote by $\xi_j$ and $\tilde \xi_j$, $j=1,\ldots, n,$ the true signal and the assimilated state. The RMSE and PCC are defined as
\begin{equation}\label{rmsepcc}
\begin{aligned}
  \mbox{RMSE} & = \sqrt{\frac{\sum_{j=1}^n(\tilde \xi_j - \xi_j)^2}{n}}, \\
   \mbox{PCC} & = \frac{\sum_{j=1}^n(\hat{\xi}_j -\overline{\hat{\xi}}_j) (\xi_j-\overline{\xi}_j)}{\sqrt{\sum_{j=1}^n(\hat{\xi}_j -\overline{\hat{\xi}}_j)^2} \sqrt{\sum_{j=1}^n(\xi_j-\overline{\xi}_j)^2}}.
\end{aligned}
\end{equation}

\subsubsection{A simple illustrative test experiment}\label{Sec:simple_example}
We start with a simple situation that contains in total $L=18$ floes in the full model with $L_0=6$ large-scale floes and $12$ small-scale floes. The small floes are parameterized as $L_s = 6$ superfloes using Algorithm \ref{alg:fnr}. The ocean contains only $8$ GB modes. Other setups are the same as those at the beginning of Section \ref{sec:uq}. The size of the ensemble in the EAKF is $1000$, which is sufficient for assimilating such a relatively simple problem without utilizing localization or covariance inflation. The observational time step is every $1.4$ hours. Three different forecast models are used in the DA:
\begin{enumerate}
  \item [1).] the full model that contains all the $18$ floes and all these $18$ floes are observed;
  \item [2).] the bare truncation model which contains only the $6$ large-scale floes and only these $6$ floes are observed; and
  \item [3).] the same truncated model and observations as (2) but including the superfloe-based noise inflation.
\end{enumerate}

\begin{figure}[h!]
\centering
	\includegraphics[width=5in]{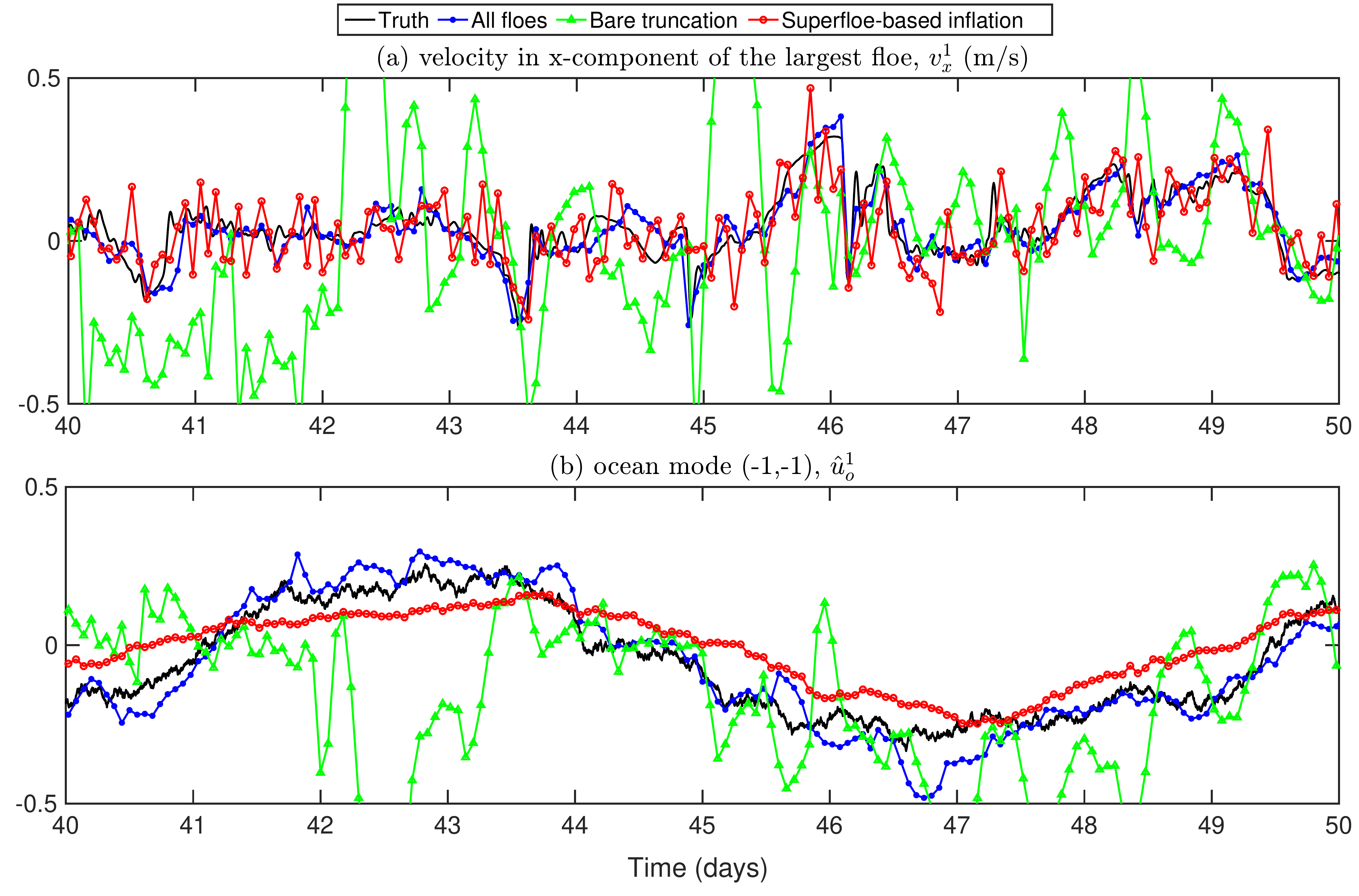}	
	\caption{The comparison of the trajectories of the assimilated velocities of the largest floe and ocean mode $(-1, -1)$. The black lines refer to the true trajectories; the blue, green, and red lines refer to assimilated trajectories when using the perfect model, bare truncation, and the superfloe-based inflation model. }
	\label{fig:f18}
\end{figure}

\begin{figure}[h!]
\centering
	\includegraphics[width=5.5in]{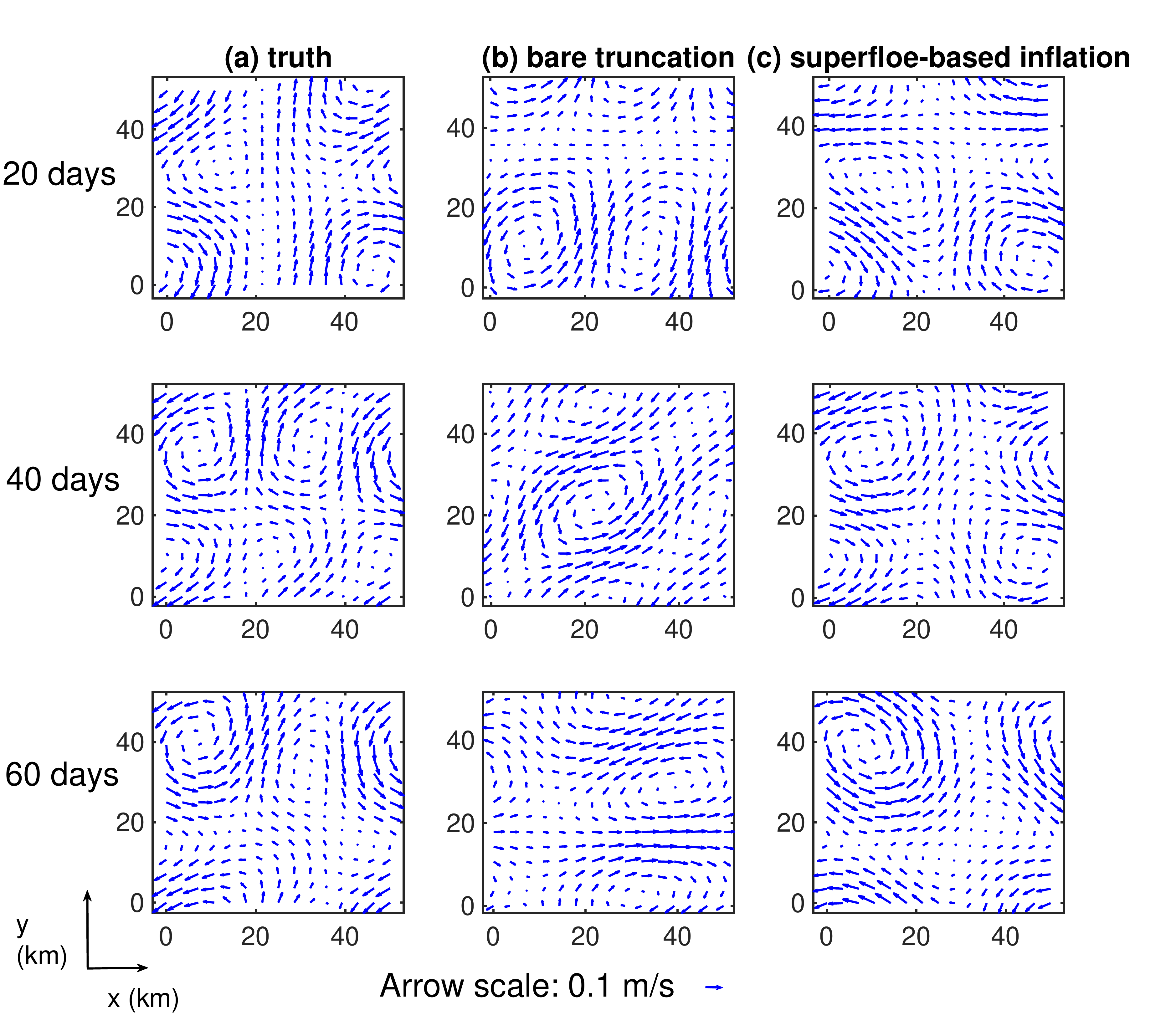}
	\caption{Snapshots of true and assimilated ocean currents when using bare truncation and superfloe-based inflation. }
	\label{fig:ocn}
\end{figure}

Figure \ref{fig:f18} shows the comparison of the trajectories of a velocity component $v^{1}_x$ and an ocean mode $(-1,-1)$. The other variables have qualitatively similar behavior. It is clear that the truncated model with noise inflation is more skillful than the one without inflation in recovering the unobserved state variables. The reason is that both the momentum and the angular momentum in \eqref{eq:pde} are driven by two things: the contact forces and the ocean drag forces. Since the bare truncation model (green curve) completely ignores the contact forces from small-scale floes, the EAKF has to treat such missing information as part of the contribution from the ocean forces. There are indeed some small fluctuations in the assimilated time series using the superfloe parameterization, which is possibly due to the noises and the inaccurate contact forces arising from the large overlaps of floes. Nevertheless, the recovered time series captures the overall time evolution of truth in an accurate fashion.  Therefore, the recovered Fourier coefficients of the ocean field contain large errors. This is also clearly indicated in the reconstructed ocean field in physical space. See Figure \ref{fig:ocn}.

\begin{figure}[h!]
\centering
	\includegraphics[width=6in]{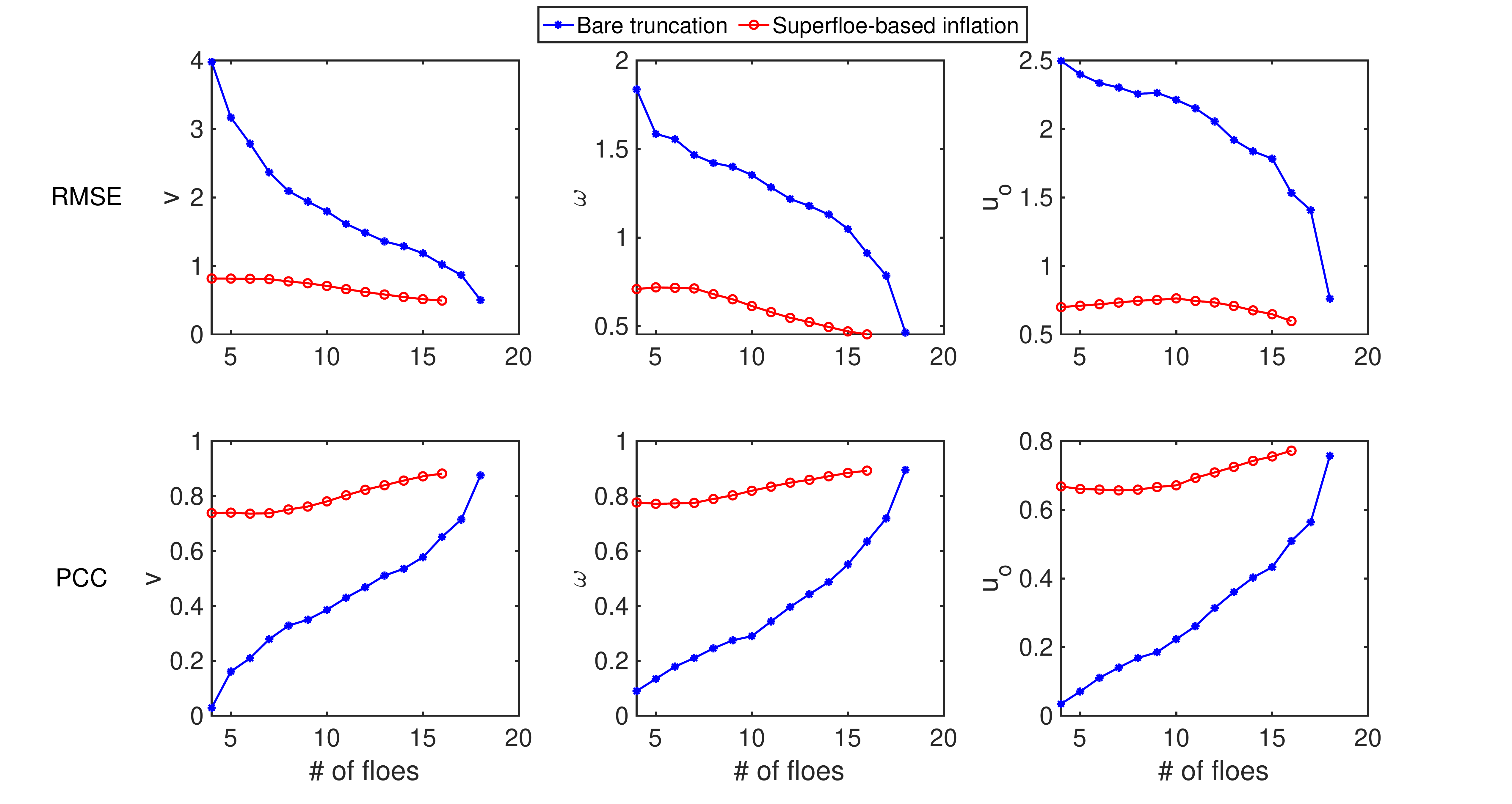}
	\caption{Comparison of the average RMSEs and PCCs of the floe velocities and ocean velocities with respect to the number of large floes.}
	\label{fig:f18ferr}
	\includegraphics[width=6in]{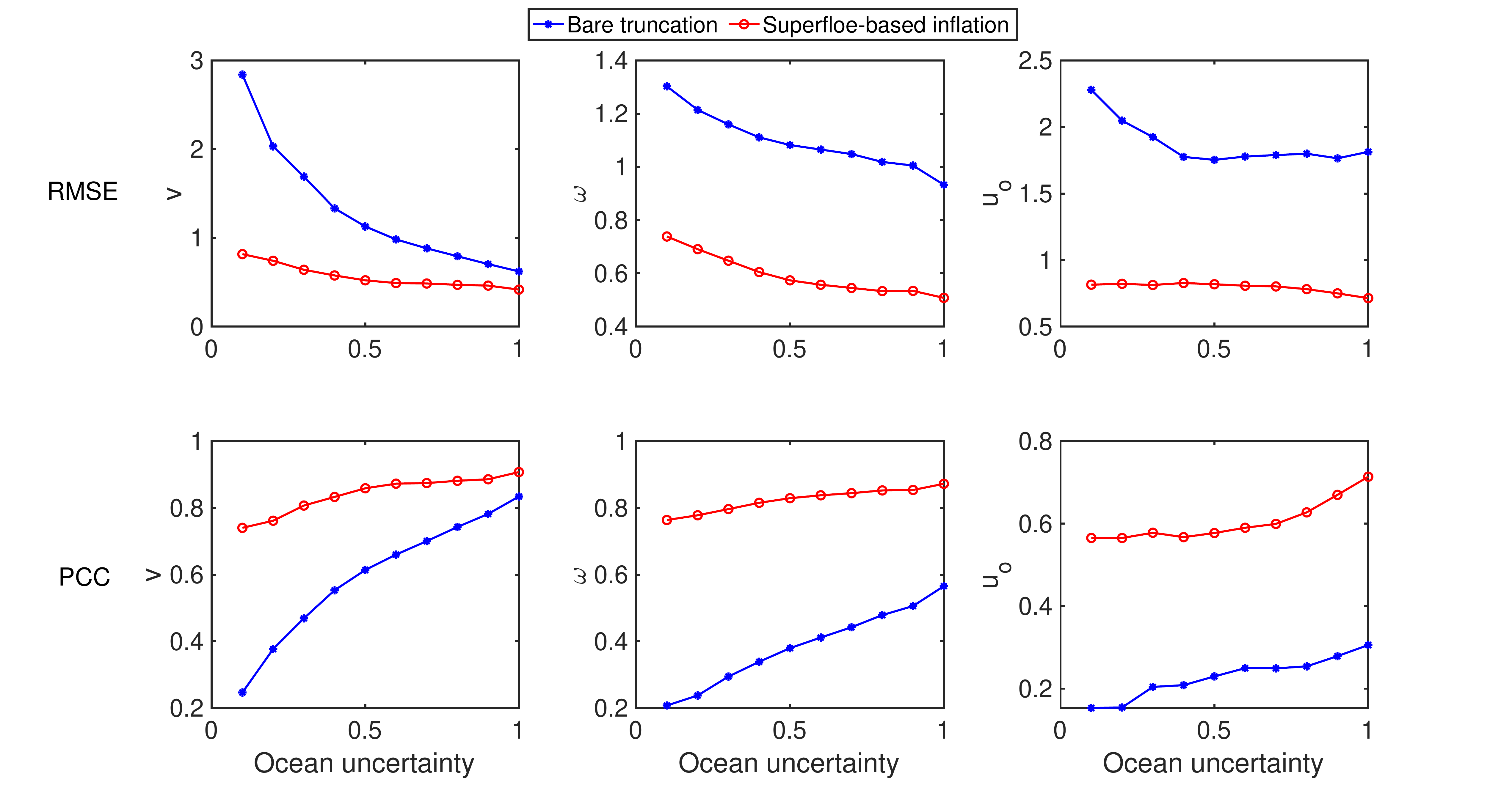}
	\caption{Comparison of the average RMSEs and PCCs of the floe velocities and ocean velocities with respect to the ocean uncertainty (noise strength).}
	\label{fig:f18oerr}
\end{figure}

Figure \ref{fig:f18ferr} shows the skill scores using the bare truncation model and the one with the superfloe-based noise inflation, where different  numbers of large-scale floes are used in the truncated models.
Here, the number of large floes varies as $L_0=4,5,\cdots,18$. Correspondingly, the number of small floes is $18-L_0$. The number of superfloes is chosen to be $\lfloor (18-L_0)/2 \rfloor$. The results in Figure \ref{fig:f18ferr} indicate that the model with the superfloe-based noise inflation consistently improves the skill scores compared with the bare truncation model. In particular, if there is only a small number of the floes retaining in the system, then the model with the superfloe-based noise inflation is significantly more skillful than the bare truncation model.

Figure \ref{fig:f18oerr} shows the skill scores with respect to the ocean uncertainty, where the noise coefficients in each ocean GB mode vary from $0.1$ to $1$.
The number of the large-scale floes and that of the superfloes are both $6$. Again, the model with the superfloe-based noise inflation outweighs the barely truncated model as a forecast model for the DA. Note that, the skill scores improve with the increase of the uncertainty of the ocean. This is because as the ocean forces increase, they dominate the contact forces, and therefore the role of the latter is weakened.

\begin{figure}[h!]
\centering
	\includegraphics[width=4in]{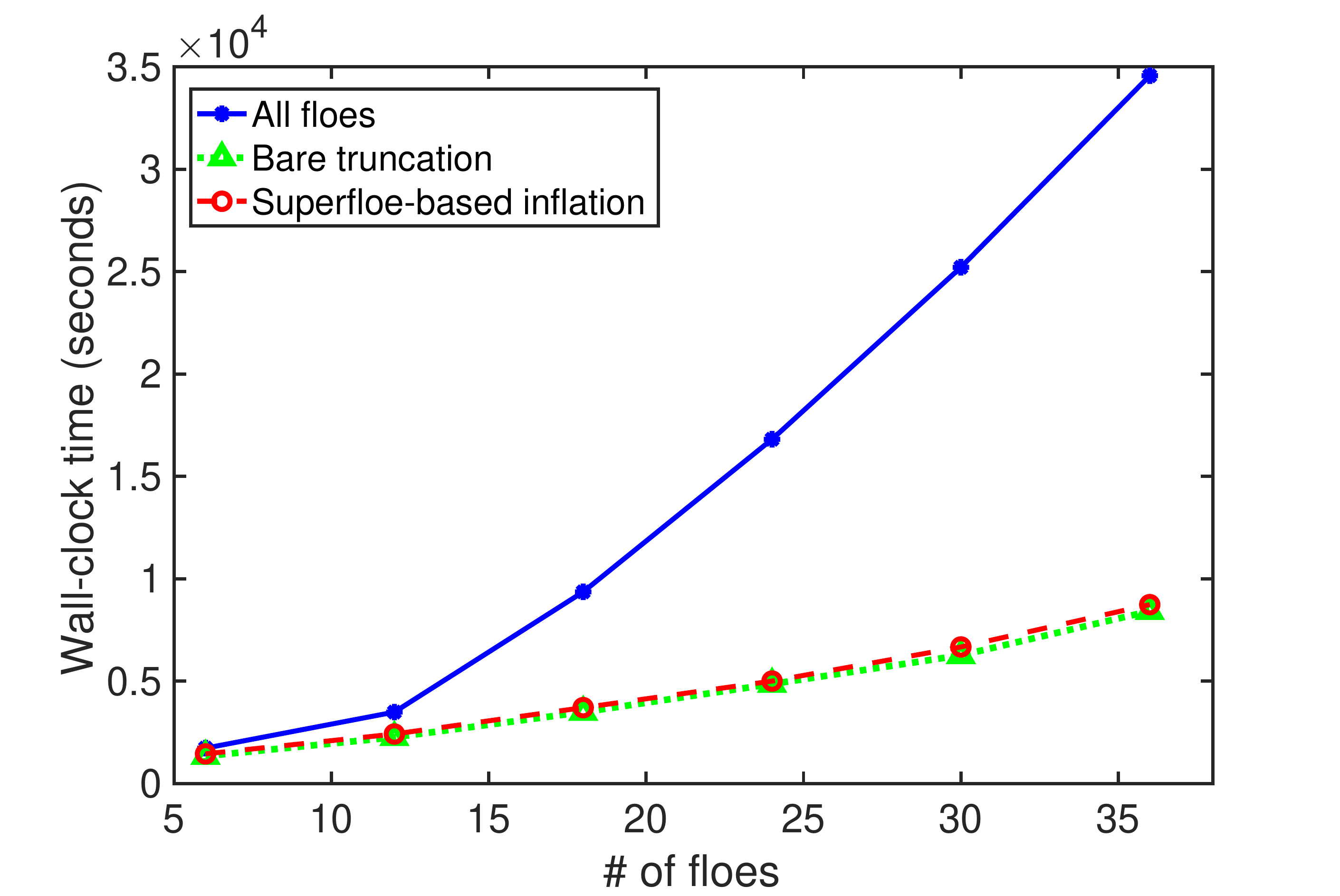}
\vspace{-0.2cm}
	\caption{Comparison of computational time cost for the case with respect to number of floes. }
	\label{fig:eakftime}
\end{figure}

To study the computational cost, which mainly depends on the number of floes in the system, we vary the number of floes as $L=6, 12, 18, 24, 30, 36$ in the full model.
The numbers of large-scale floes and the superfloes are fixed as $L_l = L_s = L/3$.
Other choices of the floe numbers are possible. Herein, we choose these numbers such that the total number of small floes are parameterized to be half number of superfloes which equals to the number of large floes in the system. 
Figure \ref{fig:eakftime} shows the corresponding comparison on the computational time.
The model with the superfloe-based noise inflation significantly reduces the computational time of DA.
Notably, since the curves of the RMSEs and PCCs in Figure \ref{fig:f18ferr} are roughly flat when the number of large floes grows, it is natural to further reduce the computational cost by keeping even fewer large-scale floes in the DA with the superfloe-based noise inflation.

\subsubsection{A more realistic test experiment}
Finally, a more realistic situation is considered. It includes more floes and a more complicated ocean field. In addition, DA with model error is taken into account here.

In this test experiment, the ocean field contains $242$ Fourier modes with $80$ GB modes and $162$ gravity modes. The noise coefficient in each GB mode is $0.1$ while that in each gravity mode is $0.02$. The values for gravity modes are smaller than those in the previous section for each Fourier mode in order to maintain the same order of the amplitude of the velocity field in the physical space. The Rossby number is still Ro$=0.1$, representing a multiscale ocean field. Since the energy of the gravity modes is relatively weak compared with that of the GB modes and the gravity modes occur in a much faster time scale, they can be treated as random perturbations on the slowly-varying GB modes. The goal here is to assimilate only the GB part of the flow, which is the typical situation in practice. Therefore, the forecast model for the Lagrangian DA excludes the gravity modes, which introduces an extra model error but accelerates the computations \cite{chen2016model}. The total number of the floes in the full system is $L=72$ with a concentration of $c=0.57$.

\begin{figure}[h!]
\centering
	\includegraphics[width=5in]{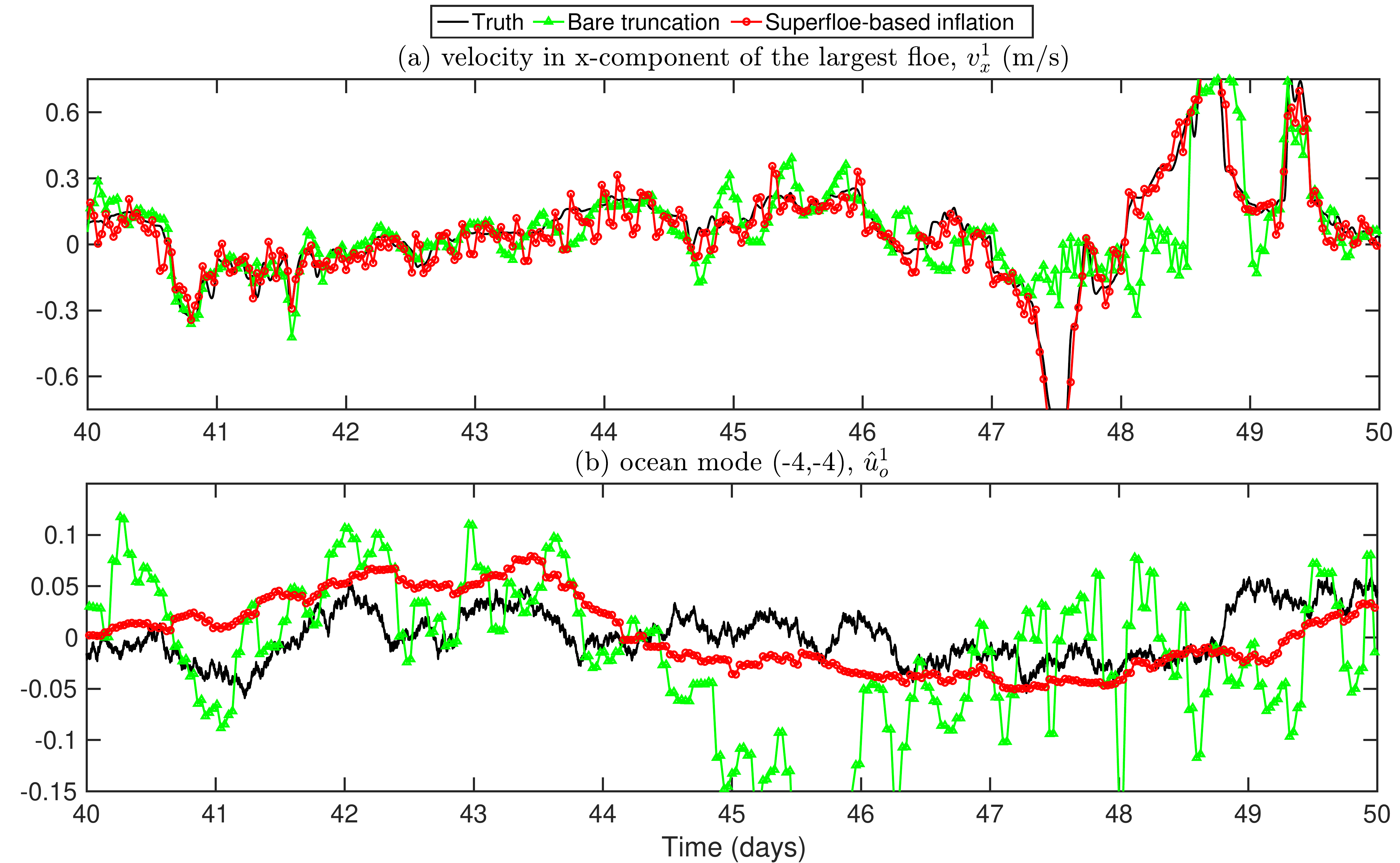}	
	\caption{The comparison of the trajectories of the assimilated velocities of the largest floe and ocean mode $(-4, -4)$. The black lines refer to the true trajectories; the blue and red lines refer to assimilated trajectories when using bare truncation and the superfloe-based inflation model. }
	\label{fig:traj}
\end{figure}

\begin{figure}[h!]
\centering
	\includegraphics[width=5.5in]{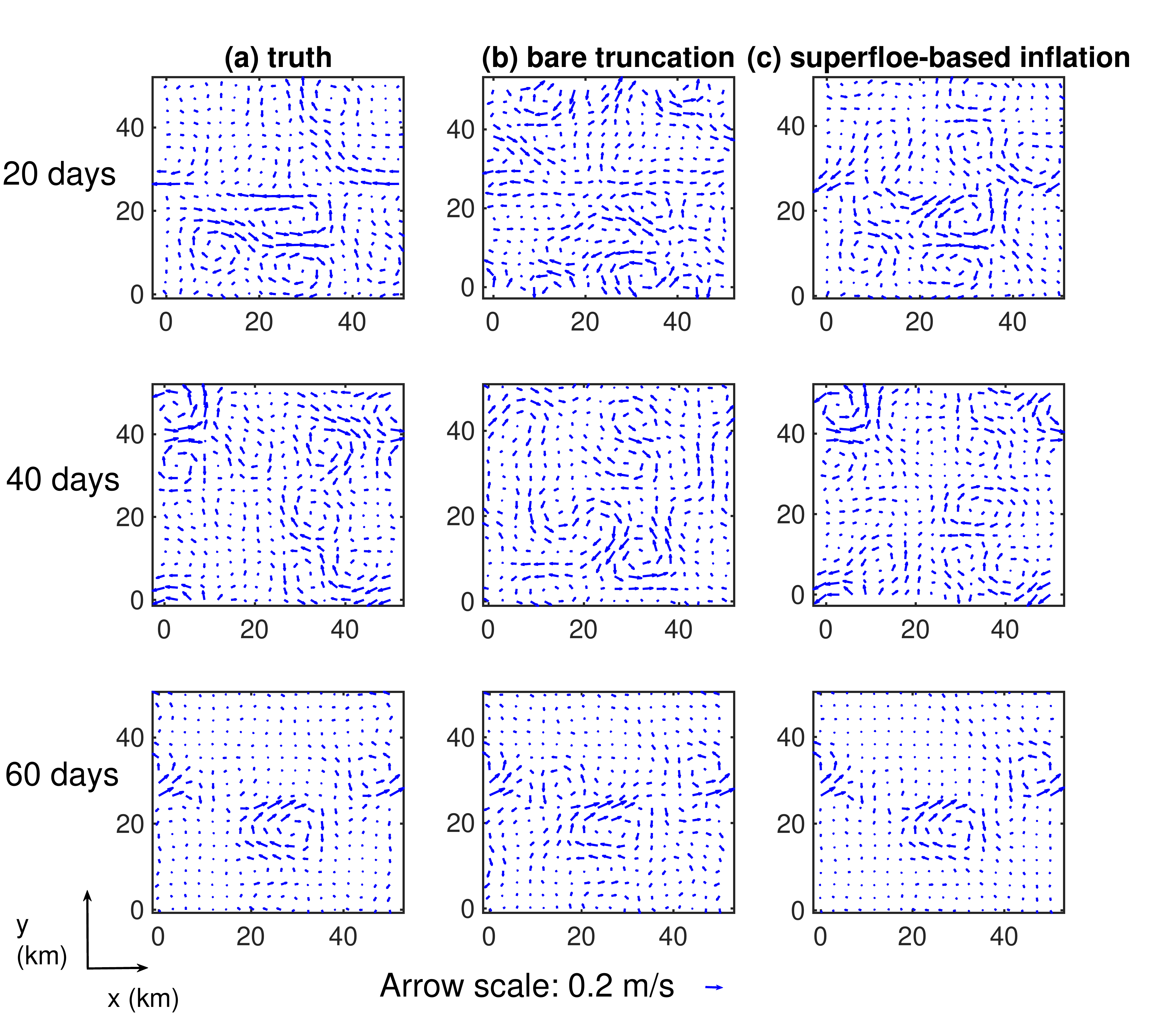}	
	\caption{Snapshots of true and assimilated ocean currents when using bare truncation and superfloe-based inflation. }
	\label{fig:ocn4}
\end{figure}

Figure \ref{fig:traj} compares the truth and the assimilated time series in terms of the velocity of the largest floe in $x$-component and the ocean mode $(-4,-4)$.
Figure \ref{fig:ocn4} shows the comparison of the truth and the reconstructed ocean flow fields. In both figures, the largest $24$ floes are retained in the reduced-order forecast model. The superfloe model, which exploits $24$ superfloes to parameterize the remaining $48$ small floes, is utilized to determine the noise inflation coefficients. Similar to the results in Section \ref{Sec:simple_example}, the bare truncation model without noise inflation leads to large errors in recovering the ocean field while the superfloe-based noise inflation significantly improves the DA skill.

\begin{figure}[h!]
\centering
	\includegraphics[width=6in]{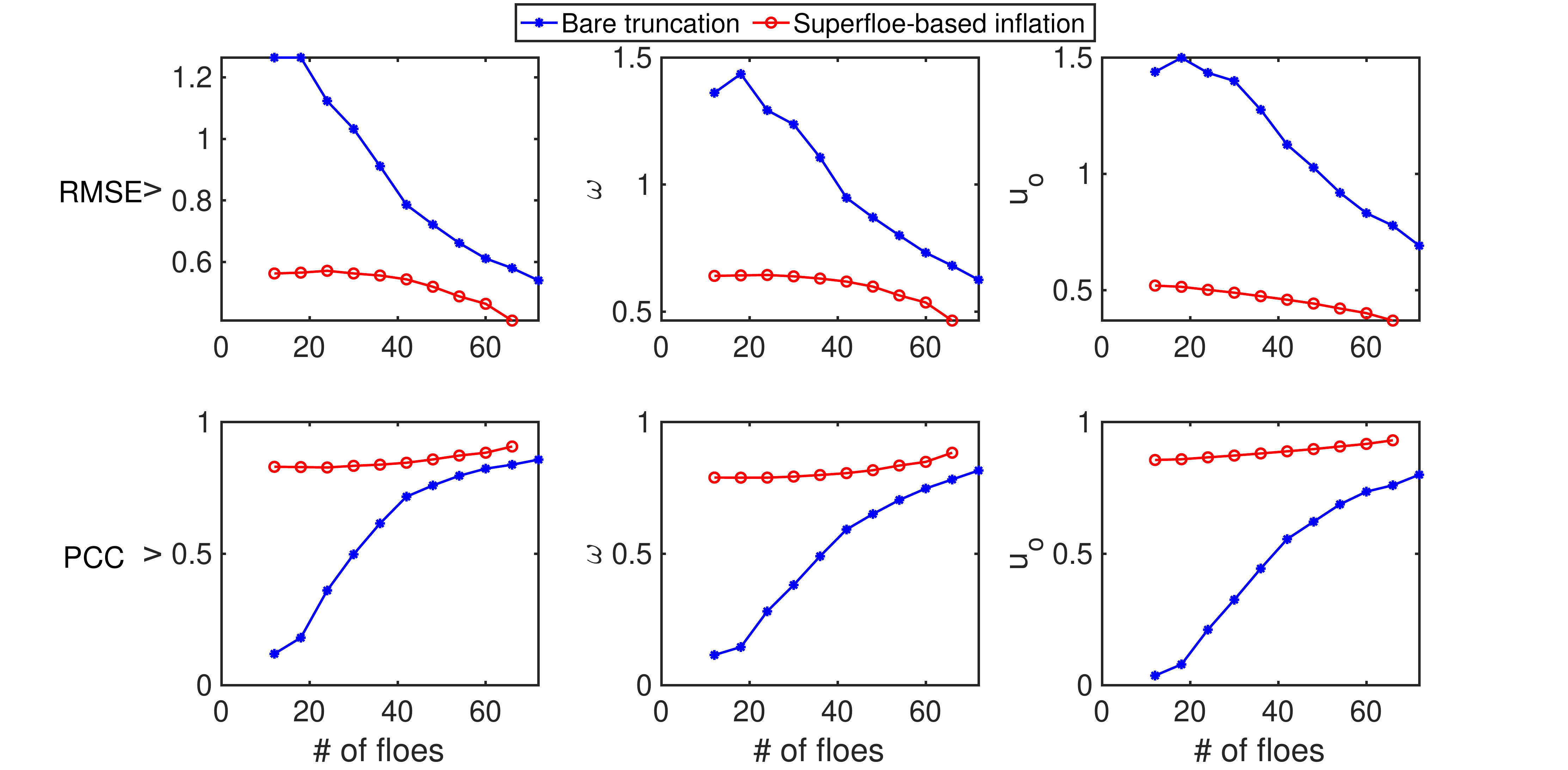}
	\caption{Comparison of the average RMSEs and PCCs of the floe velocities and ocean modes with respect to the number of large floes.}
	\label{fig:err}
\end{figure}

Figure \ref{fig:err} shows the skill scores as a function of the number of large-scale floes in the reduced-order system. The number of the superfloes is set to be half as many as the number of the small-scale floes that are unresolved in the reduced-order model. The error associated with the bare truncation model without noise inflation increases dramatically when the number of the large-scale floes decreases while the error associated with the model using the superfloe-based noise inflation remains at a low level. One interesting finding is that applying the reduced-order model with superfloe-based noise inflation using only a small number of the large-scale floes even outweighs the one that includes all the $72$ floes. In fact, the forecast model here does not include the gravity modes. Nevertheless, the superfloe-based noise inflation automatically takes into account such an effect. Thus, the noise inflation compensates both the contact forces from the small-scale floes and the model error due to the ignorance of the gravity modes.

\section{Concluding remarks} \label{sec:con}
In this paper, an efficient physics-constrained superfloe parameterization is developed that significantly reduces the computational cost of the DEM model for sea ice.
The superfloe parameterized system captures the main features of sea ice floe dynamics as well as the long-term non-Gaussian statistical features. It also facilitates a systematic noise inflation scheme that advances the ensemble-based DA algorithms.
Yet, there remain a few simplifications in the current setups of the superfloe paramterization. One assumption is that the groups of small floes would remain bundled together indefinitely. The associated validity can potentially depend on the initial distribution of floes and distribution of surrounding sea forces. One natural direction of improving the framework here is to systematically develop a criterion for breaking large pieces of the superfloes into small ones. This should also be incorporated with the melting and the fracture features, which are currently not included in the DEM utilized here.
Other future work includes applying the superfloe parameterization to more realistic sea ice models with more features and their associated DA problems.
To capture more physics of the Arctic sea ice dynamics,
a sophisticated model would include more features such as fracturing and the impact of the salinity, temperature, and atmosphere.
For these more complicated setups, the EAKF should be combined with the localization and covariance inflation for effectively assimilating the coupled system. These techniques can be naturally incorporated into the current framework.
Finally, although it has been shown the results of the superfloe parameterization with different numbers of large-scale floes and superfloes, developing systematic strategies to determine these numbers with certain theoretic guidelines for UQ and DA remains as future work.  

\section*{Acknowledgement}
N.C. and S.N.S are partially funded by ONR MURI N00014-19-1-2421. Q.D. is a Van Vleck Visiting Assistant Professor under this grant.


\bibliographystyle{siam}
\bibliography{ref}

\appendix{}


\section{Model parameters}
Table \ref{tab:param} includes the parameters and their physical units of the idealized DEM model \eqref{eq:pde}. These parameters correspond to the physical variability in the marginal ice zone in the Arctic ocean.

\begin{table}[ht]
\centering
\begin{tabular}{c|cc}
\hline
& Simulation domain & periodic, 50 km$\times50$ km \\
& Numerical scheme & Euler-Maruyama \\
& Time-marching step size & $\Delta t = 25$ seconds \\
Simulation  & Simulation final time & $T \sim 120$ days \\
& Number of DA ensembles & 1000 \\
& Observation variables & $\textbf{x}, \Omega$ \\
& Observation noise strength & $\sigma_\textbf{x}=$80m, $\sigma_\Omega = $0.01rad \\
& Observational time-step size & $100\Delta t \approx 1.4$ hours \\ \hline

& Sea ice density & $\rho_o = 900$kg/m$^3$ \\
& Size (radius $r$) distribution & $p(r) =\frac{1.5 }{r^2}$ \\
& Thickness ($h$) distribution & $p(h)= 0.59 h e^{-0.77 h}$ m \\
& Size typical range & $r \in [1 \text{km}, 10 \text{km}]$ \\
Floe & Thickness typical range & $h \in [0.1 \text{m}, 3.5 \text{m}]$ \\
& Concentration typical range & $[0.1, 0.8]$ \\
& Shear and Young's modulus & $E^{lj}=G^{lj}=1.25\times10^8$ Pa \\
& Coulomb friction & $\mu^{lj} = 0.2$ \\
& Velocity scale & $\sim 0.1$m/s \\
& Angular velocity scale & $\sim 10^{-5}$rad/s \\ \hline

& Seawater density & $\rho_o = 10^3$kg/m$^3$ \\
& Velocity scale  & $U_o\sim0.1$m/s \\
& Ocean drag coefficient &$d_o = 3\times10^{-3}$ \\
Ocean& Rossby number & Ro$=0.1$ \\
& Damping coefficients for the ocean modes & $d_{\mathbf{k},\zeta} = \phi_{\mathbf{k},\zeta}= 0.5$ \\
& Long-term mean forcing for the GB modes & $f_{\mathbf{k},{\zeta}}=0.1\exp(i2\pi/14t)$  \\
& Long-term mean forcing for the gravity modes & $f_{\mathbf{k},{\zeta}}=0$  \\
& Ocean GB mode uncertainty strength & $\sigma_{\mathbf{k},\zeta}=0.1$ \\
& Ocean gravity mode uncertainty strength & $\sigma_{\mathbf{k},\zeta}=0.05$ \\  \hline
 \end{tabular}

 \vspace{0.5cm}
\caption{Parameters and their physical units of the idealized DEM model \eqref{eq:pde}.}
\label{tab:param}
\end{table}

\end{document}